\documentclass[12pt]{article}

\usepackage[T1]{fontenc}
\usepackage{amsmath}
\usepackage{amssymb}
\usepackage{textcomp}
\input{diagrams}
\diagramstyle{h=21pt,w=21pt,scriptlabels,midshaft,PostScript=dvips,nohug,shortfall=4pt}
\usepackage[noBBpl]{mathpazo}
\usepackage{graphicx}

\renewcommand{\texttt}[1]{{\fontfamily{pcr}\fontseries{m}\fontshape{n}\selectfont#1}}

\usepackage[applemac]{inputenc}

\renewcommand{\ldots}{\relax\ifmmode\ldotp\ldotp\ldotp\else$\m@th\ldotp\ldotp\ldotp\ $\fi}
\providecommand{\qed}{\hspace*{\fill}\nolinebreak[1]\hspace*{\fill} $\Box$}
\providecommand{\nomoreqed}{\renewcommand{\qed}{}}

\renewcommand{\epsilon}{\varepsilon}

\newcommand{\into}			{\hookrightarrow}
\newcommand{\isopil}{\stackrel{\raisebox{0.1ex}[0ex][0ex]{\(\sim\)}}%
			{\raisebox{-0.15ex}[0.28ex]{\(\rightarrow\)}}}

\newcommand{\upperstar}{^{\raisebox{-0.25ex}[0ex][0ex]{\(\ast\)}}}
\newcommand{\lowerstar}{_{\raisebox{-0.33ex}[-0.5ex][0ex]{\(\ast\)}}}

\newcommand{\df}{\: {\raisebox{0.255ex}{\normalfont\scriptsize :\!\!}}=}

\newcommand{\tensor}	{\otimes}

\newcommand{\Hom}	{\operatorname{Hom}}
\newcommand{\End}	{\operatorname{End}}

\newcommand{\Id}	{\operatorname{Id}}
\newcommand{\ov}{\overline}
\newcommand{\fat}[1]{\mathbf{{#1}}}

\newcommand{\overskrift}[1]{\par\noindent\relax{\LARGE #1}\par\bigskip}

\newcounter{dummycounter}

\newenvironment{punkt-a}%
{%
	\begin{list}%
	{(\alph{dummycounter})\hfill}%
	{\usecounter{dummycounter}%
	\setlength{\itemsep}{0em}\setlength{\parsep}{0em}\setlength{\topsep}{0em}%
	\setlength{\itemindent}{0em}\setlength{\labelwidth}{1.5em}
	\setlength{\labelsep}{0.3em}\setlength{\leftmargin}{1.8em}}%
}%
{\end{list}}

\newenvironment{punkt-i}%
{%
	\begin{list}%
	{(\roman{dummycounter})}%
	{\usecounter{dummycounter}%
	\setlength{\itemsep}{0em}\setlength{\parsep}{0em}\setlength{\topsep}{0em}%
	\setlength{\itemindent}{0em}\setlength{\labelwidth}{1.8em}%
	\setlength{\labelsep}{0.6em}\setlength{\leftmargin}{2.4em}}%
}%
{\end{list}}

\newcommand{\hovedfont}{\normalfont\bfseries}
\usepackage{theorem}
\theorembodyfont{\normalfont\itshape}
\theoremheaderfont{\hovedfont}
	\theoremstyle{change}
\newtheorem{lemma}{Lemma.}[section]
\newtheorem{prop}[lemma]{Proposition.}
\newtheorem{satz}[lemma]{Theorem.}
\newtheorem{cor}[lemma]{Corollary.}
	\theorembodyfont{\normalfont}

\newtheorem{BM}[lemma]{Remark.}

\newtheorem{taller}[lemma]{$\!\!$}
\newenvironment{blanko}[1]%
{\begin{taller}{\hovedfont #1}\normalfont}%
{\end{taller}}
{%
\begin{list}{\em Definition. }%
{\setlength{\labelsep}{0mm}\setlength{\leftmargin}{0mm}%
\setlength{\labelwidth}{0mm}\setlength{\listparindent}{\parindent}%
\setlength{\parsep}{\parskip}\setlength{\partopsep}{0mm}}%
\item%
}%
{%
\end{list}%
}
\newenvironment{dem}%
{%
\begin{list}{\em Proof. }%
{\setlength{\labelsep}{0mm}\setlength{\leftmargin}{0mm}%
\setlength{\labelwidth}{0mm}\setlength{\listparindent}{\parindent}%
\setlength{\parsep}{\parskip}\setlength{\partopsep}{0mm}}%
\item%
}%
{%
\qed\end{list}%
}
\newenvironment{dem*}[1]%
{%
\begin{list}{\em #1 }%
{\setlength{\labelsep}{0mm}\setlength{\leftmargin}{0mm}%
\setlength{\labelwidth}{0mm}\setlength{\listparindent}{\parindent}%
\setlength{\parsep}{\parskip}\setlength{\partopsep}{0mm}}%
\item%
}%
{%
\qed\end{list}%
}
{%
\begin{list}{\em Proof. }%
{\setlength{\labelsep}{0mm}\setlength{\leftmargin}{0mm}%
\setlength{\labelwidth}{0mm}\setlength{\listparindent}{\parindent}%
\setlength{\parsep}{\parskip}\setlength{\partopsep}{0mm}}%
\item%
}%
{%
\qed\end{list}%
}
\newenvironment{bevis*}[1]%
{%
\begin{list}{\em #1 }%
{\setlength{\labelsep}{0mm}\setlength{\leftmargin}{0mm}%
\setlength{\labelwidth}{0mm}\setlength{\listparindent}{\parindent}%
\setlength{\parsep}{\parskip}\setlength{\partopsep}{0mm}}%
\item%
}%
{%
\qed\end{list}%
}

\newenvironment{blanko*}[1]%
{%
\begin{list}{\bf {#1} }%
{\setlength{\labelsep}{0mm}\setlength{\leftmargin}{0mm}%
\setlength{\labelwidth}{0mm}\setlength{\listparindent}{\parindent}%
\setlength{\parsep}{\parskip}\setlength{\partopsep}{0mm}}%
\item%
}%
{%
\end{list}%
}

\usepackage{mathrsfs}

\providecommand{\lastUpdate}[1]{#1}

\newcommand{\dotsum}{\begin{picture}(15,5)(-7.5,-2.5)
\put(0,-2){\makebox(0,0){{\scriptsize \(\bullet\)}}}
\put(0,1.5){\makebox(0,0){{\(+\)}}}
\end{picture}}
\newcommand{\smalldotsum}{\begin{picture}(9,3)(-5,-1.5)
\put(0,-2){\makebox(0,0){{\tiny \(\bullet\)}}}
\put(0,1.5){\makebox(0,0){{\scriptsize\(+\)}}}
\end{picture}}

\providecommand{\kat}[1]{\text{\textbf{\textsl{#1}}}}

\newcommand{\topile}{\raisebox{-1.5pt}{\(\stackrel{\rTo}{\rTo}\)}}
\newcommand{\topileback}{\raisebox{-1.5pt}{\(\stackrel{\lTo}{\lTo}\)}}

\newcommand{\op}{^{\text{{\rm{op}}}}}
\newcommand{\Ob}{\operatorname{Ob}}
\newcommand{\Arr}{\operatorname{Arr}}

\newcommand{\isleftadjointto}{\dashv}

\providecommand{\B}{\mathbb{B}}
\providecommand{\S}{\mathbb{S}}
\providecommand{\T}{\mathbb{T}}

\providecommand{\Z}{\mathbb{Z}}

\newcommand{\CC}{\mathscr{C}}
\newcommand{\DD}{\mathscr{D}}

\renewcommand{\SS}{\mathscr{S}}

\newcommand{\Grpd}{\kat{Grpd}}

\newcommand{\Set}{\kat{Set}}

\newcommand{\Top}{\kat{Top}}

\newcommand{\sSet}{\kat{sSet}}
\newcommand{\Cat}{\kat{Cat}}

\newcommand{\sCat}{\kat{sCat}}
\newcommand{\CCat}{\kat{CCat}}
\newcommand{\Ho}{\text{{\rm{Ho}}}}
\newcommand{\semiCat}{\kat{\textonehalf Cat}}

\newcommand{\nCat}{\kat{nCat}}
\newcommand{\disju}{{\textstyle{\coprod}}}

\newcommand{\Dmono}{\Delta_{\operatorname{mono}}}

\newcommand{\smallcoprod}[2]{\overset{#2}{\underset{#1}{\textstyle{\coprod}}}} 
\newarrow{Lig}{=}{=}{=}{=}{=}

\newarrow{Mono}{>}{-}{-}{-}{>}
\newarrow{Epi}{-}{-}{-}{-}{>>}

\usepackage{texdraw}
\input{txdtools}

\everytexdraw{%
	\drawdim pt 
	\textref h:C v:C
}

\newcommand{\Bdot}{%
 \drawdim pt 
 \fcir f:0 r:2
}

\newcommand{\hop}{%
 \drawdim pt 
 \rmove (0 10)
}

\newcommand{\korthop}{%
 \drawdim pt 
 \rmove (0 7.5)
}

\newcommand{\Blink}{%
 \drawdim pt 
 \linewd 1.5
 \savecurrpos (*ex *ey)
 \rmove (0 2.2)
 \rlvec (0 5.6)
 \move (*ex *ey)
 \rmove (0 10)
}
\newcommand{\kortBlink}{%
 \drawdim pt 
 \linewd 1.5
 \savecurrpos (*ex *ey)
 \rmove (0 2.2)
 \rlvec (0 3.1)
 \move (*ex *ey)
 \rmove (0 7.5)
}

\setbox2=\hbox{%
\raisebox{0.5pt}{%
\begin{texdraw} \drawdim pt 
	\setunitscale 1
	\move (0 0)
	\Bdot
\end{texdraw}}}
\newcommand{\intextB}{\,\usebox2\,}
\setbox3=\hbox{%
\raisebox{-1.5pt}{%
\begin{texdraw} \drawdim pt 
	\setunitscale 1
	\move (0 0)
	\Bdot\korthop\Bdot
\end{texdraw}}}
\newcommand{\intextBB}{\,\usebox3\,}

\setbox6=\hbox{%
\raisebox{-1.5pt}{%
\begin{texdraw} \drawdim pt 
	\setunitscale 1
	\move (0 0)
	\Bdot\kortBlink\Bdot
\end{texdraw}}}
\newcommand{\intextBbB}{\,\usebox6\,}

\newcommand{\intextBBB}{%
  \raisebox{-3.5pt}{
\begin{texdraw} \drawdim pt 
	\setunitscale 1
	\move (0 0)
	\Bdot\korthop\Bdot\korthop\Bdot
\end{texdraw}
}
}

\newcommand{\intextBbBbB}{%
  \raisebox{-3.5pt}{
\begin{texdraw} \drawdim pt 
	\setunitscale 1
	\move (0 0)
	\Bdot\kortBlink\Bdot\kortBlink\Bdot
\end{texdraw}
}
}

\def\scaleFactor{10} % MUST BE AN INTEGER

\newcommand{\dropglob}[1]{%
\setlength{\unitlength}{0.003\DiagramCellWidth}
\multiply \unitlength by \scaleFactor
\begin{picture}(0,0)(0,0)
\qbezier(-28,-4)(0,-18)(28,-4)
\put(0,-14){\makebox(0,0)[t]{$\scriptstyle {#1}$}}
\put(28.6,-3.7){\vector(2,1){0}}
\end{picture}
}

\newcommand{\topglob}[1]{%
\setlength{\unitlength}{0.003\DiagramCellWidth}
\multiply \unitlength by \scaleFactor
\begin{picture}(0,0)(0,0)
\qbezier(-28,11)(0,25)(28,11)
\put(0,21){\makebox(0,0)[b]{$\scriptstyle {#1}$}}
\put(28.6,10.7){\vector(2,-1){0}}
\end{picture}
}

\newcommand{\lift}[2]{%
\setlength{\unitlength}{1pt}
\begin{picture}(0,0)(0,0)
\put(0,{#1}){\makebox(0,0)[b]{${#2}$}}
\end{picture}
}

\newcommand{\cel}[1]{\ensuremath{\mathsf{#1}}}

\newarrow{Lig}{=}{=}{=}{=}{=}

\newcommand{\To}{\begin{diagram}[w=0.5em]{}&\rTo&{}\end{diagram}}

\makeatletter
\renewcommand{\ps@headings}
  {\setlength{\headheight}{41pt}%
   \setlength{\headsep}{12pt}%
   \renewcommand{\@oddhead}{\parbox{\textwidth}{%
      \scriptsize
      \texttt{\lastUpdate {2005-11-28}
      \hfill Joachim Kock: Weak identity arrows in higher categories 
      \hfill  [\thepage/\pageref{lastpage}]}
      \\ \rule[8pt]{\textwidth}{0.3pt}}%
   }
  \renewcommand{\@oddfoot}{}
  \renewcommand{\@evenfoot}{}%
}
\makeatother

\begin{document}

% Assuming that the package graphicsx has already been loaded
% HERE IS THE SYMBOL FOR FAT DELTA:
% this must come after begin docment:
\setbox1=\hbox{\includegraphics{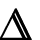}}
\providecommand{\grosdelta}{\usebox1}
\providecommand{\gros}{\grosdelta}
\providecommand{\grosop}{\grosdelta{\!}\op}
% WARNING: THIS METHOD IS INCOMPATIBLE WITH USING \alignat*
% See example in the file factorisation.tex
\renewcommand{\SS}{\kat{S}}

\newcommand{\ord}[1]{\mathbf{#1}}
\newcommand{\nCob}{\kat{nCob}}
\newcommand{\SFairCat}{\SS\kat{-FairCat}}
\newcommand{\SCat}{\SS\kat{-Cat}}

\pagestyle{headings}

\vspace*{24pt}
\begin{center}
  
\overskrift{Weak identity arrows in %\\[4pt]
higher categories}

\bigskip

\noindent
\textsc{Joachim Kock}
\bigskip

% \footnotesize Universit\'e de Nice Sophia-Antipolis
% \normalsize

\end{center}

\begin{abstract}
  There are a dozen definitions of weak higher categories, all of
  which loosen the notion of composition of arrows.  A new approach is
  presented here, where instead the notion of identity arrow is
  weakened --- these are tentatively called fair categories.  The
  approach is simplicial in spirit, but the usual simplicial category
  $\Delta$ is replaced by a certain `fat' delta of `coloured
  ordinals', where the degeneracy maps are only up to homotopy. 
  The first part of this exposition is aimed at a broad
  mathematical readership and contains also a brief introduction to
  simplicial viewpoints on higher categories in general.  It is
  explained how the definition of fair $n$-category is almost forced
  upon us by three standard ideas. 
  
  The second part states some basic results about fair categories, and
  give examples, including Moore path spaces and cobordism categories.
  The category of fair $2$-categories is shown to be
  equivalent to the category of bicategories with strict composition
  laws.  Fair $3$-categories correspond to tri\-categories with strict
  composition laws.  The main motivation for the theory is Simpson's
  weak-unit conjecture according to which $n$-groupoids with strict
  composition laws and weak units should model all homotopy $n$-types.
  A proof of a version of this conjecture in dimension $3$ is
  announced, obtained in joint work with A.~Joyal.  Technical details
  and a fuller treatment of the applications will appear elsewhere.
\end{abstract}

\bigskip

\setcounter{tocdepth}{1}

\makeatletter

\renewcommand{\tableofcontents}{%
   \begin{center}
\begin{minipage}{124mm}
   \begin{center}
		\bf{\contentsname}
	\end{center}
   \footnotesize
   \begin{center}
		\@starttoc{toc}
	\end{center}	
\end{minipage}
	\end{center}
	\addvspace{3em \@plus\p@}
}

\renewcommand{\section}{\@startsection {section}{1}{\z@}%
{-3.5ex \@plus -1ex \@minus -.2ex}%
{2.3ex \@plus.2ex}%
{\normalfont\large\bfseries}}

\renewcommand*{\l@section}[2]{%
  \ifnum \c@tocdepth >\z@
    \addpenalty\@secpenalty
%     \addvspace{0.5em \@plus\p@}%
    \setlength\@tempdima{1.5em}%
    \begingroup
      \parindent \z@ \rightskip \@pnumwidth
      \parfillskip -\@pnumwidth
      \leavevmode %\bfseries
      \advance\leftskip\@tempdima
      \hskip -\leftskip
      #1\nobreak\hfil \nobreak\hb@xt@\@pnumwidth{\hss #2}\par
    \endgroup
  \fi}

  \renewcommand{\subsection}{\@startsection {subsection}{1}{\z@}%
{-3.5ex \@plus -1ex \@minus -.2ex}%
{2.3ex \@plus.2ex}%
{\normalfont\normalsize\bfseries}}

\makeatother

\tableofcontents

\pagebreak

\begin{flushright}
  \begin{minipage}{3.6in}
    \em How would you rate the choice of terminology?
    
    $\Box$ Excellent \quad
    $\Box$  Good \quad
    $\boxtimes$  Fair \quad
    $\Box$ Weak \quad
    $\Box$ Poor
  \end{minipage}
\end{flushright}

\setcounter{section}{-1}
%%%%%%%%%%%%%%%%%%%%%%%%%%%%%%%%%%%%%%%%%%%%%%%%%%
\section{Introduction}
%%%%%%%%%%%%%%%%%%%%%%%%%%%%%%%%%%%%%%%%%%%%%%%%%%

\begin{blanko*}{Higher categories.}
  While conventional category theory encompasses mathematics in the
  paradigm of objects and arrows-between-objects, higher category
  theory considers also $2$-arrows between arrows, $3$-arrows between
  $2$-arrows, and so on.  Higher categorical structures were first
  discovered in algebraic geometry and homotopy theory, and are now
  becoming increasingly important in many areas of mathematics, as
  well as in theoretical physics and computer science.  It is not
  difficult to define strict higher categories, but the crucial point for
  applicability of higher categories is {\em weakening}, and the
  theory is inherently of homotopical nature.  The theory of weak
  higher categories is young, with about a dozen competing definitions
  (cf.~Leinster's survey~\cite{Leinster:ten}).  All these
  approaches emphasise weakened composition laws.
\end{blanko*}

\begin{blanko*}{Weak identity arrows.}
  The present paper introduces a new approach to the problematic, with
  a systematic theory of weak identity arrows.  
  While strict identity arrows appear naturally whenever the arrows 
  represent some sort of mappings, there are higher-dimensional
  contexts where they seem less nature-given and sometimes problematic:
  
  --- identity arrows are degenerate things, and tend to collapse
  other things as well (cf.~the Eckmann-Hilton argument);
  
  --- identity arrows tend to be of a `wrong geometrical type' (for
  instance: not cofibrant, or of defective dimension (e.g., identity
  $n$-cobordisms are `cylinders' of height zero)).
  
  One would like instead to speak of merely `up-to-homotopy' identity
  arrows, so some sort of homotopical context is needed, i.e.~a 
  category with a notion of equivalence $\simeq$.

  \bigskip
    
  In a usual category the identity arrows appear as the
  image of a degeneracy map
  \begin{diagram}[w=5ex,h=3.5ex,tight]
    O & \rTo  & A
  \end{diagram}
  where $O$ is the set of objects, and $A$ the space of arrows.  (And
  this map is subject to the well-known identity axioms relating
  to the composition map.)
  
  There are essentially two possible ways one can try to define weak
  identity arrows.  Either one can stick with the map $O \to A$, but
  weaken the axioms.  This amounts to making non-canonical choices for
  the identity arrows, and handling the ensuing coherence issues.  Such
  approaches to weak higher categories are usually called `algebraic',
  since the operations (in this case the `nullary' operation $O \to
  A$) remain well-defined.
  
  The other approach, which is the route taken in this work, is to weaken
  the very shape of the diagram above: the idea is to avoid the 
  artificial choices and instead relax the structure
  by specifying an acyclic space of weak identity arrows $U$, sitting
  between $O$ and $A$ like this:
  \begin{diagram}[w=5ex,h=3.5ex,tight]
  O &   & A  \\
  \uTo<\simeq  & \ruTo    &   \\
  U  &   & 
  \end{diagram}
  (and satisfying certain axioms relative to the composition law).
  Thus the theory remains purely diagrammatic, and becomes homotopical
  rather than algebraic, in the sense that there is no longer a
  well-defined `nullary' operation --- it has been replaced by an
  up-to-homotopy-only operation.  Given the space $U$, it is sometimes
  possible to revert to an algebraic description by merely choosing a
  pseudo-section $u:O\to U$; the axioms satisfied by $U$ then
  automatically equips $u$ with the extra structure needed to serve as
  weak identity arrows and with coherence constraints built in.
  Conversely, given some algebraic notion of weak identity arrows,
  ideally their existence should be a property, not a structure, in an
  up-to-homotopy sense, meaning that the space $U(x)$ of all possible
  weak identity structures on an object $x$ should be contractible.
  
  One example to have in mind (treated in detail in
  Section~\ref{Sec:examples}) is the following.
  % Let $X$ be a topological space, and let $O$ denote its set of points.
  A Moore path in a topological space $X$ is a 
  continuous map $I\to X$ where $I$ is an
  interval of any positive length.  
  % Let $A$ denote the space of all
  % Moore paths; we have source and target maps $O \topileback A$. 
  Moore paths can be composed in a strictly associative way (by
  concatenation of intervals, hence without need of passing to
  homotopy classes of paths), but since we exclude the interval of
  zero length there are no units for this composition.  Composing a
  path with its `inverse' suggests that the weak units should be the
  null-homotopic loops, but we need to specify in which way each such
  null-homotopic loop is a weak unit: a weak unit at a point $x$
  is a null-homotopic loop at $x$ {\em together} with a null-homotopy.
  The space $U(x)$ of such pairs is contractible, and our $U$ is the
  disjoint union $\coprod U(x)$ where $x$ runs over the points of $X$.
\end{blanko*}

\begin{blanko*}{Tamsamani categories and fair categories.}
  The diagrammatic weakening just outlined is analogous to the
  weakening of the composition laws in the theory of Tamsamani higher 
  categories \cite{Tamsamani:thesis}, which seems to be the most 
  developed theory of higher categories.  While the strict composition
  law is a diagram of shape
  \begin{diagram}[w=5ex,h=3.5ex,tight]
  A & \lTo & A \times_O A ,
  \end{diagram}
  in a Tamsamani category this shape has been replaced by
    \begin{diagram}[w=5ex,h=3.5ex,tight]
  A &   & A\times_O A  \\
    & \luTo    & \uTo>\simeq  \\
    &   & B
  \end{diagram}

  Tamsamani's theory is iterated simplicial, and the weakening of the
  composition laws is obtained by relaxing the nerve condition (see
  \ref{Tams} below).  It is not possible to obtain a similar
  diagrammatic
  weakening of the identity arrow structure in a purely simplicial
  setting, because the identity structure does not come from an
  external condition like the nerve condition; it is inherent from
  the degeneracies in $\Delta$.
  Instead, the usual simplicial category $\Delta$ is replaced
  by a larger category $\gros$ which is introduced in this paper.
  Still the theory is simplicial in spirit, and it goes hand in hand
  with Tamsamani's theory: one version of it combines the two
  weakenings.

  In order to emphasise the weak-identity-arrow aspects of
  the theory, we will mostly concentrate on the case of strict
  composition laws.
  Such categories (strict composition laws but weak
  identity arrows) are tentatively called {\em fair categories}.
  Strict composition laws are not usually what is encountered in
  nature, but according to a strong form of Simpson's conjecture, the
  homotopy category of weak $n$-categories is equivalent to the
  homotopy category of $n$-categories with strict composition laws and
  only weak identity arrows.  A proof of a weak version of Simpson's
  conjecture in dimension $3$ is announced here.  This result has been
  obtained jointly with Andr\'e Joyal and will appear separately
  \cite{Joyal-Kock:traintracks}.  The emphasis on strict composition
  laws is justified by this result.  Section~\ref{Sec:groupoids} 
  contains some background for Simpson's conjecture and motivation for
  the theory.
\end{blanko*}

\begin{blanko*}{The fat delta.}
  The fat delta $\gros$, which is the basis for the theory, has the following 
concise description as a subcategory of the category of arrows in 
$\Delta$: its objects are the epimorphisms in $\Delta$, and the
arrows are the monomorphisms in $\Arr(\Delta)$;  these are the 
commutative diagrams in $\Delta$ whose vertical arrows are epimorphisms
and whose top arrow is a monomorphism:
\begin{diagram}[w=3.6ex,h=3.6ex,tight]
\cdot & \rInto  & \cdot  \\
\dEpi  &    & \dEpi  \\
\cdot  & \rTo  & \cdot
\end{diagram}

  However, this description does not convey much intuition.  Instead,
  a less concise but more conceptual definition is given, and it is
  shown how the fat delta and the definition of fair category follow
  quite naturally from three standard ideas: the simplicial idea, the
  homotopy idea, and the idea of `semi-'.  The fat delta is the
  category of coloured semi-ordinals.
\end{blanko*}

%%%%%%%%%%%%%%%%%%%%%%%%%%%%%%%%%%%%%%%%%%%%%%%%%%
\subsection*{Organisation of the paper and overview of the results}
%%%%%%%%%%%%%%%%%%%%%%%%%%%%%%%%%%%%%%%%%%%%%%%%%%

The first three sections are aimed at a broad mathematical readership
and serve to introduce the simplicial viewpoint on higher categories,
establish some terminology, and to motivate and introduce the fat
delta.  With these preparations, the definition of fair category
\ref{def} is a one-liner: {\em A fair category in $\SS$ is a functor
$\grosop\to\SS$ that preserves equimorphisms, discrete objects, and
fibre product over discrete objects.}

The second part of the paper, containing the results,
is a bit more advanced, and some of the
details of proofs are deferred to forthcoming papers.
  
Section for section:

\begin{blanko*}{\S 1 \ Nerves and simplicial enrichment.}
  The definition of fair categories is simplicial in spirit, and the
  first section recalls the basic viewpoint of categories as special
  simplicial sets. $n$-categories are defined inductively as
  certain special simplicial objects in the category of
  $(n-1)$-categories.  The fundamental combinatorial structure is
  $\Delta$, the category of non-empty finite ordinals.
\end{blanko*}

\begin{blanko*}{\S 2 \ Homotopy and coloured categories.} \
  In order to be able to talk of weakness and `up-to-homotopy', we
  must work in a category $\SS$ with a notion of equivalence --- these are
  called {\em coloured categories} in this work.  Weakening the
  composition laws in the simplicial viewpoint leads to
  Tamsamani's notion of higher category, cf.~\cite{Tamsamani:thesis}.
%   In order to carry along the notion of equivalence in the induction 
%   step, some technical arguments on its relation to {\em discrete objects}
%   are needed.  These are given in the appendix.
  
  Those two ideas: `simplicial' and `up-to-homotopy', lead directly to
  the concept of {\em coloured ordinals}.  The
  category of (finite, nonempty) coloured ordinals plays the same
  r\^ole for coloured categories as the usual $\Delta$ does for plain
  categories.
\end{blanko*}

%   This means that
%   certain `coloured' morphisms in $\SS$ (called {\em equimorphisms})
%   are singled out, and the game is to treat them as if they were
%   invertible.  Typical examples are $\Top$, the category of
%   topological spaces and continuous maps, colouring the weak homotopy
%   equivalences, and $\Cat$, the category of (small) categories and
%   functors, colouring the `equivalences' (i.e., full, faithful, and
%   essentially surjective functors).  (For the theory to work, $\SS$ is
%   furthermore required to have a notion of discrete objects, similar
%   to the meaning of this term in $\Top$ and $\Cat$.)

\begin{blanko*}{\S 3 \ Semi-categories, coloured semi-categories,
  and the `fat delta' $\gros$.}
  To arrive at fair categories, a third idea is involved, namely {\em
  semi-categories} (categories without identities).  Combining this
  idea with the previous two ideas yields the (coloured) category
  $\gros$ of (finite, nonempty) coloured semi-ordinals.  This
  category captures the shape of fair categories, in the sense that
  fair categories will be certain $\grosop$-diagrams in $\SS$, just as
  usual categories are seen as $\Delta\op$-diagrams.
\end{blanko*}

\begin{blanko*}{\S 4 \ Definition of fair categories and fair $n$-categories.}
  A {\em fair category} in $\SS$ is by definition a functor $X:\grosop
  \to \SS$ that preserves equimorphisms, discrete objects, and fibre
  products over discrete objects.  That $X$ preserves equimorphisms
  expresses the up-to-homotopy identity arrow condition.  Preserving
  fibre products over discrete objects expresses the condition of
  strict composition laws.  (At this point the reader may wish to have
  a glance at the examples given in Section~\ref{Sec:examples}.)
  
  Inductively, a fair $n$-category is obtained with
  $\SS=\kat{(n-1)Cat}$, the category of fair $(n-1)$-categories.  
  For this to make sense we must define what an equivalence of fair 
  $\SS$-categories is.  Some
  technicalities on discrete objects needed in the definition are
  relegated to an appendix --- the reader can safely think of the
  usual notion of discrete objects in $\Top$ or in $\Cat$. 
  
  One basic result, also needed in the
  induction, is that {\em every fair $1$-category is isomorphic to a
  usual category.}
\end{blanko*}

\begin{blanko*}{\S 5 \ Fair $2$-categories.}
  The following is the main result of \S 5:
\begin{blanko*}{Proposition \ref{B}.} {\em
  The category of fair $2$-categories is equivalent to the category
  of bicategories with strict composition laws.}
\end{blanko*}
  The fair $2$-category viewpoint on such a 
  bicategory encodes all the possible unit structures on the 
  underlying semi-bicategory.  
\end{blanko*}

\begin{blanko*}{\S 6 \ Fair $3$-categories.}
  Morally, fair $3$-categories should correspond to tri\-categories
  with strict composition laws.  A special semi-strict case of
  particular interest is worked out: 
  \begin{blanko*}{Proposition~\ref{corr}.} {\em
    Fair monoidal strict $2$-categories correspond to monoidal 
    strict $2$-categories with weak units in the sense of
    Gordon-Power-Street~\cite{Gordon-Power-Street}.}
  \end{blanko*}
\end{blanko*}

\begin{blanko*}{\S 7 \ $n$-groupoids, homotopy $n$-types, and Simpson's conjecture.}
  This section serves as a second introduction, explaining the main
  motivation for considering weak units, and the justification for
  considering only strict composition laws:
  Simpson~\cite{Simpson:9810} has conjectured that every homotopy
  $n$-type arises as the geometric realisation of a strict
  $n$-groupoid with weak identity arrows --- for a suitable notion of
  weak identity arrows.  There is a straightforward notion of fair
  $n$-groupoid, and it is conjectured that this notion will fulfil
  Simpson's conjecture.  A proof of a version of the conjecture in
  dimension $3$ is announced here, obtained jointly with Andr\'e
  Joyal.  The key result is this:
  \begin{blanko*}{Theorem~\ref{JK}.} {\em
    (Joyal-Kock~\cite{Joyal-Kock:traintracks}.)
    Every braided monoidal category arises as $\End(I)$, where $I$ 
    is a weak unit in an otherwise completely strict monoidal $2$-category.}
  \end{blanko*}
  \begin{blanko*}{Corollary~\ref{JKcor}.} {\em
    (Cf.~\cite{Joyal-Kock:traintracks}.) Strict $2$-groupoids with invertible tensor product and weak units
    can model all $1$-connected homotopy $3$-types.}
  \end{blanko*}
\end{blanko*}

\begin{blanko*}{\S 8 \ Examples.}
  Weak identity arrows arise as `honest' replacements for
  `artificial' identity arrows or approximation to non-existent 
  identity arrows.  Three examples of this are given, the first in 
  some detail, the other two more succinctly.
  
  The first example concerns Moore paths in a topological space.  The
  weak identity arrows are null-homotopic paths together with a
  null-homotopy.
    
  The second example concerns a monoidal model category in the sense
  of Hovey~\cite{Hovey:model}.  If the unit is cofibrant then the full
  subcategory of all cofibrant objects is a genuine monoidal coloured
  category.  If the unit is not cofibrant one gets instead a fair
  monoidal coloured category, the space of weak units being the
  cofibrant replacements of the unit.  
%   (This example is the idea behind a
%   proof of a non-linear version of Deligne's Hochschild cohomology
%   conjecture \cite{Kock-Toen:0304} (joint with To\"en), although in
%   the end we actually found a simpler way.)

  The final example concerns cobordism categories, and is perhaps
  of particular interest.  In summary:
  \begin{blanko*}{Proposition~\ref{CobProp}.} {\em
    Oriented $n$-cobordisms naturally assemble into a fair Tamsamani
    $2$-category, for which the straight cylinders are weak identity 
    arrows.}
  \end{blanko*}
\end{blanko*}

\begin{blanko*}{Acknowledgements.}
  The bulk of this paper was written during the Summer of 2003, and
  preliminary versions of it were circulated in connection with three
  talks: in Aarhus in April 2003, at the CATS2 conference in Nice, May
  2003, and at the Workshop on Higher-order Geometry and
  Categorification in Lisbon, July 2003.  In Lisbon, I presented an
  example of a fair $2$-monoid in $\Cat$ and claimed that it was
  equivalent to the braid category.  This statement is wrong --- I am
  grateful to Andr\'e Enriques for pointing out the problems in my
  argument.  This led to the discovery that all fair $2$-fold monoidal
  categories collapse to symmetric ones \cite{Kock:commutativity}, and
  for a long time it seemed that Simpson's conjecture would be false;
  for this reason the present paper was shelved.  Decisive new insight
  resulted from many conversations with Andr\'e Joyal during my year
  in Montréal, culminating with Theorem~\ref{JK}.

  This research was supported by a Marie Curie Fellowship 2001--2002, by the
  University of Nice (2003), and by a CIRGET postdoc grant at the University
  of Qu\'ebec at Montr\'eal, 2004.  I am very much indebted to Andr\'e
  Hirschowitz, Bertrand To\"en, and Andr\'e Joyal for their help and
  encouragement, and I would like also to acknowledge fruitful
  conversations and e-mail correspondence with Anders Kock, Bill Lawvere, 
  Carlos Simpson, Clemens Berger, Markus Spitzweck, and Tom Leinster.
\end{blanko*}

% \pagebreak

%%%%%%%%%%%%%%%%%%%%%%%%%%%%%%%%%%%%%%%%%%%%%%%%%%
\section{Categories as simplicial sets}
%%%%%%%%%%%%%%%%%%%%%%%%%%%%%%%%%%%%%%%%%%%%%%%%%%

We shall briefly review Tamsamani's definition of higher categories,
and settle on some basic terminology.  The fundamental viewpoint is 
that a category is a simplicial set satisfying certain conditions.

\begin{blanko}{Simplicial sets, categories, and the (strict) Segal 
  condition.}
  Let $\Delta$ be the category whose objects are
  the nonempty finite ordinals, (i.e.~the linearly ordered nonempty
  finite sets $$ \ord 0 = \{0\}, \quad \ord 1 = \{0 \leq 1\}, \quad
  \ord 2 = \{0\leq 1\leq 2\}, \quad \ldots ,\quad \ord n = \{0\leq
  1\leq 2\leq \ldots\leq n\},\quad\ldots)
$$
% \begin{eqnarray*}
%   \ord 0 & = & \{0\},   \\
%   \ord 1 & = & \{0 \leq 1\}, \\
%   \ord 2 & = & \{0\leq 1\leq 2\}, \\
%   \vdots & & \\
%   \ord n & = & \{0\leq 1\leq 2\leq \ldots\leq n\}, \\
%   \vdots & &
% \end{eqnarray*}
and whose arrows are the order-preserving maps, i.e., functions 
$f : \ord m\to \ord n$ such that
$f(i)\leq f(j)$ in $\ord n$ whenever $i\leq j$ in $\ord m$.
It is convenient to interpret each $\ord n$ as a category by
viewing each inequality $i \leq j$ as an arrow from $i$ to $j$. 
Then $\Delta$ becomes the full subcategory of $\Cat$ formed by the
categories $\ord 0, \ord 1, \ord 2, \ldots$

A simplicial set is a functor $\Delta\op\to\Set$; a simplicial map is 
a natural transformation of such functors.  We denote by $\sSet$ the 
category of simplicial sets and maps.
  
The {\em nerve} of a (small) category $\mathbf{C}$ is by definition the 
simplicial set
% the set of objects and the set of all
% arrows,
%   $$
%   C_0 \df \Ob(C) \qquad 
%     C_1 \df \coprod_{x,y\in C_0} \Hom_C(x,y) 
%   $$
%   fit naturally into a simplicial set, called the {\em nerve} of $C$:
%   \begin{diagram}[w=7ex,tight,shortfall=2ex,scriptlabels]
%     C_0 & \pile{\lTo^t \\ \rTo[shortfall=3ex] \\ \lTo_s}  & C_1
%     &  \pile{\lTo \\ \rTo[shortfall=3ex] \\ \lTo \\\rTo[shortfall=3ex] \\ \lTo}  &
%     C_2  & \dots  
%   \end{diagram}
%   The maps $s,t$ associate 
%   to a given arrow its source and target respectively.  The map $C_0 
%   \to C_1$ associates to each object its identity 
%   arrow.  The set $C_2$ is the set of composable arrows:
%   $$
%   C_2 \df C_1 \times_{C_0} C_1 = \coprod_{x,y,z\in C_0} \big(\Hom_C(x,y)
%   \times \Hom_C(y,z)\big) ,
%   $$
%   and the middle map $C_2 \to C_1$ is the composition law, etc.
%   
%   The whole construction is captured in the
%   the nerve of a category $C$
% is simply
\begin{eqnarray*}
C:	\Delta\op & \longrightarrow & \Set  \\
	{}\ord n & \longmapsto & \Hom_{\Cat} (\ord n, \mathbf{C}) .
\end{eqnarray*}
There are natural identifications
  $$
  C_0 = \Ob(\mathbf{C}) , \qquad 
    C_1 = \coprod_{x,y\in C_0} \Hom_{\mathbf{C}}(x,y) ,
  $$
  and more generally $C_k$ is interpreted as the set of strings of 
  $k$ composable arrows in $\mathbf{C}$.
  Functors between
  categories turn into simplicial maps between their nerves, and the
  whole construction defines a functor from categories to simplicial
  sets.  This functor is fully
  faithful.
  The simplicial sets that arise as nerves of categories are 
  characterised by {\em strict Segal condition}: the natural maps
%   \begin{equation}\label{strictSegal}
%       C_k \simeq C_1 \times_{C_0}\dots\times_{C_0} C_1,
%   \end{equation}
  \begin{equation}\label{strictSegal}
      C_{p+q} 
      \to C_p \times_{C_0} C_q, \quad\forall p,q,
  \end{equation}
  are isomorphisms.
  (Here the fibre product refers to the maps $C_p \to C_0\leftarrow 
  C_q$, `target of last arrow' and `source of first arrow',
  respectively.)
%   (Equivalently,
%   $C_{p+q} \simeq C_p \times_{C_0} C_q, \forall p,q$.) 
 Hence it makes
  sense to say that {\em a category is a simplicial set satisfying the
  strict Segal condition}.

  (This viewpoint and the terminology {\em nerve} go back to
  Grothendieck~\cite{Grothendieck:exp190}, 1959.  
  The
  (non-strict) Segal condition played a fundamental r\^ole in
  Segal~\cite{Segal:cat-coh} (1974), and was named after
  him by Tamsamani~\cite{Tamsamani:thesis} in 1996.)

It will be helpful to have the following graphical description
in mind: The category $\ord n$ is pictured with a dot for each
object;   the arrows are not drawn explicitly; instead
the order is expressed by having the dots arranged in a column, with
$0$ at the bottom and $n$ at the top.  Now a functor $\ord m \to \ord
n$ (i.e., an order preserving map) is represented by linking each dot
in $\ord m$ to its image dot in $\ord n$; the order preservation then
corresponds to the requirement that these strands do not cross over
each other.  Here is a picture of the most fundamental arrows in
$\Delta$:

\label{delta-graphic}
\begin{center}\label{4}
\begin{texdraw}
	\setunitscale 1.1
	\move (0 0)
	% THIS IS IDENTITY
	\bsegment
	\move (0 -5) \Bdot\hop\Bdot
	\move (20 0) \Bdot
	\linewd 0.2 \move (0 -5) \lvec (20 0) \lvec (0 5)
	\htext (10 -18){$\ord 1 \to \ord 0$}
	\esegment
	\move (70 0)
	% THIS IS SOURCE
	\bsegment
	\move (0 0) \Bdot
	\move (20 -5) \Bdot\hop\Bdot
	\linewd 0.2 \move (0 0) \lvec (20 -5)
	\htext (10 -18){$\ord 0 \to \ord 1$}
	\esegment
	\move (140 0)
	% THIS IS TARGET
	\bsegment
	\move (0 0) \Bdot
	\move (20 -5) \Bdot\hop\Bdot
	\linewd 0.2 \move (0 0) \lvec (20 5)
	\htext (10 -18){$\ord 0 \to \ord 1$}
	\esegment
	\move (210 0)
	% THIS IS COMPOSITION
	\bsegment
	\move (0 0) \Bdot\hop\Bdot
	\move (20 -5) \Bdot\hop\Bdot\hop\Bdot
	\linewd 0.2 \move (0 0) \lvec (20 -5)
		    \move (0 10) \lvec (20 15)
	\htext (10 -18){$\ord 1 \to \ord 2$}
	\esegment
\end{texdraw}
\end{center}
If $C:\Delta\op\to \Set$ is (the nerve of) a category, then
the four arrows in the picture are mapped to the
following four maps: the map $C_0 
\to C_1$ that sends an object to its identity arrow; the source map $C_1\to C_0$;
the target map $C_1 \to C_0$; and the 
composition of arrows $C_2 \to C_1$.

The Segal condition is an exactness condition: In $\Delta$, two arrows
$\ord m \leftarrow \ord 0 \to \ord n$ admit a pushout if and only if the dot of
$\ord 0$ is included in $\ord m$ as the last dot, and in $\ord n$ as the first
dot (or conversely), like in this example:
  \begin{center}
  \begin{texdraw}
%       \arrowheadsize l:4 w:3 \arrowheadtype t:V
	  \setunitscale 0.8
	  \move (0 0) \bsegment
	  \move (0 60) \Bdot
	  \lvec (5 5) \Bdot \lvec (70 -10) \Bdot
	  \move (0 60) \lvec (65 45) \Bdot 
	  \lvec(70 -10)
	  
	  \move (-5 -5) \Bdot \lvec (60 -20) \Bdot
	  
	  \move (75 55) \Bdot \lvec (80 0) \Bdot
	  \move (85 65) \Bdot \lvec (90 10) \Bdot
	  \move (95 75) \Bdot \lvec (100 20) \Bdot
	  \esegment
	  \htext (200  30) {\begin{diagram}[w=32pt,h=24pt,tight]
	  {}\ord 0 & \rTo  & \ord 3  \\
	  \dTo  &    & \dTo  \\
	  {}\ord 1  & \rTo  & \NWpbk \relax \ord 4
	  \end{diagram}
	  }
	  \move (240 0)
  \end{texdraw}
  \end{center}
  These amalgamated sums in $\Delta$ are fibre products in
  $\Delta\op$, and the Segal condition amounts precisely to saying
  that $X:\Delta\op\to\Set$ preserves fibre products over $\ord 0$
  (i.e., those fibre products over $\ord 0$ that happen to exist).
\end{blanko}

\begin{blanko}{Strict higher categories.}\label{strict-higher}
  A (strict) $2$-category is defined to be a category enriched over
  $\Cat$.  This means that the definition of category is repeated, but
  the hom sets (and the structure maps between them) are replaced by
  hom categories $\Hom_C(x,y)$ (and the structure maps by functors).
  There is a category $\kat{2Cat}$ whose objects are $2$-categories,
  and whose arrows are $2$-functors.  Inductively then, a (strict)
  $n$-category is a category enriched over $\kat{(n-1)Cat}$.
  Unwinding this recursive definition, one finds that an $n$-category
  has objects, arrows between objects, $2$-arrows between arrows, and
  so on up to $n$-arrows between $(n-1)$-arrows. $k$-arrows
  can be composed in $k$ compatible ways, provided their
  lower-dimensional cells match appropriately, and for each $k$-arrow
  ($k<n$) there is an identity $(k+1)$-arrow.
  
  It is important to note that throughout this induction, the objects
  of an $n$-category always form a {\em set}, not an $(n-1)$-category.
  Now a set can in a natural way be considered as a category, namely a
  category whose only arrows are the identity arrows.  Similarly a
  category can be considered as a $2$-category, etc., and in
  particular a set can be considered as an $n$-category for any $n$;
  these are just the $n$-categories whose only $k$-arrows are the 
  identity arrows ($0<k\leq n$).
  Such $n$-categories are called {\em discrete}.  With this remark,
  the iterated enrichment can be conveniently formulated in terms of
  nerves, and this formulation is the key to the generalisation and
  weakening we will come to in the next section: A {\em $0$-category}
  is just a set, and a $1$-category is just a usual category.  An {\em
  $n$-category} is a functor
  $$
  C : \Delta\op\to \kat{(n-1)Cat}
  $$
  such that $C_0$ is discrete and such that the strict Segal
  condition~(\ref{strictSegal}) holds.
\end{blanko}

\begin{blanko}{$\SS$-categories.}\label{SCat}
  The definition of higher categories involves considering simplicial
  objects in many different categories and imposing the `category
  condition'.  Hence the natural need for an abstract notion of
  $\SS$-category: this should be a simplicial object
  $X:\Delta\op\to\SS$ in a category $\SS$, such that $X_0$ is discrete
  and such that $X_{p+q} \simeq X_p \times_{X_0} X_q$.  For this to
  make sense, $\SS$ must be have a notion of {\em discrete objects},
  and it should possess {\em fibre products over discrete objects}.
  These two notions are unavoidable, and in fact the category
  condition is nothing but preservation of these two notions: 
  Namely, define $\ord 0$ to be the only discrete object in $\Delta$ (or
  in $\Delta\op$).  (Thinking of discreteness as dual to
  connectedness, clearly this is the only reasonable way to define
  discrete objects in $\Delta$.)

  {\sc Definition.} Let $\SS$ be a category with a notion of discrete
  objects, and admitting fibre products over discrete objects.  An
  {\em $\SS$-category} is a functor $X:\Delta\op\to\SS$ that preserves
  discrete objects and fibre products over discrete objects.  Let
  $\SCat$ denote the category whose objects are $\SS$-categories $X:
  \Delta\op\to\SS$ and whose morphisms are the natural transformations
  between such.
\end{blanko}

\begin{blanko}{Weakening.}
  An important lesson from category theory is summarised in the
  following two remarks.  On level 0: when comparing two sets (i.e,
  $0$-categories), it is more important to say if they are isomorphic
  than if they are actually equal.  On level 1: given two categories
  it is more important to say if they are equivalent than if they are
  actually isomorphic.  Recall that an equivalence of categories is a
  functor that is fully faithful and essentially surjective.
  Henceforth such a functor will rather be called an {\em
  equimorphism}.  According to this lesson, the above definition of
  $2$-category is not the `correct' one, since the strict Segal
  condition refers to isomorphism of categories.  A weaker notion of
  higher category is obtained by requiring the Segal maps $C_{p+q} \to
  C_p \times_{C_0} C_q$ to be merely equimorphisms, not isomorphisms:
\end{blanko}

\begin{blanko}{Tamsamani higher categories.}\label{Tams}
  Define a Tamsamani $0$-category to be a set.
  A {\em weak $n$-category in the sense of 
  Tamsamani}~\cite{Tamsamani:thesis}
  is defined inductively as a functor
  $$
  C:\Delta\op\to \kat{(n-1)wCat}
  $$
  such that $C_0$ is discrete, and satisfying the (non-strict) Segal
  condition, namely that the morphisms $C_{p+q} \to
  C_p \times_{C_0} C_q$ should be equimorphisms in $\kat{(n-1)wCat}$, the 
  category of weak $(n-1)$-categories.  `Equimorphism' means
  `fully faithful' and `essentially surjective', notions which are
  also defined inductively.  These definitions rely crucially on 
  properties of the notion of discrete object, which we come to in a 
  minute.
  
  Note that an equimorphism is not in general 
  invertible, so there is no longer any well-defined composition 
  like this:
  \begin{diagram}[w=6ex,h=4ex,tight]
    C_1  & \lTo & C_1\times_{C_0} C_1
  \end{diagram}
  This map is now defined only `up to homotopy': it exists
  only inasmuch as we regard the equimorphism $C_2 \isopil C_1 
  \times_{C_0} C_1$ as invertible.  The new structure is rather
  this:
  \begin{diagram}[w=6ex,h=4ex,tight]
    C_1 &   & C_1\times_{C_0} C_1 \\
    & \luTo    & \uTo>\simeq  \\
    &   & C_2
  \end{diagram}

  The theory of Tamsamani higher categories seems to be the most
  developed among the theories of higher categories, thanks in
  particular to the work of Simpson and his collaborators (see for
  example \cite{Simpson:internalhom}, \cite{Simpson:9810058},
  \cite{Hirschowitz-Simpson}, \cite{Pellissier}, 
  \cite{Toen-Vezzosi:0212}).  The main reason is that it is a
  simplicial theory, and huge bodies of simplicial methods in
  homotopical algebra can be applied.
\end{blanko}

%%%%%%%%%%%%%%%%%%%%%%%%%%%%%%%%%%%%%%%%%%%%%%%%%%
\section{Homotopy and coloured categories}
%%%%%%%%%%%%%%%%%%%%%%%%%%%%%%%%%%%%%%%%%%%%%%%%%%

\begin{blanko}{Coloured categories.}\label{CCat}
  Weakening makes sense in categories equipped with
  a notion of equivalence.  Then a property can be said to hold up
  to equivalence (or up to homotopy), meaning that it only holds
  inasmuch as equivalences are considered as equalities.  The
  following terminology is handy, but not standard in the literature.  A
  {\em coloured category} is a category $C$ with a specified
  subcategory $W$ comprising all the objects.  The arrows in $W$ are
  called {\em coloured arrows} or {\em equimorphisms} or {\em equiarrows}.
  An {\em
  equivalence} is a zigzag of equimorphisms, and two objects in $C$
  are {\em equivalent} if there is an equivalence between them.  (This
  relation is clearly reflexive, symmetric, and transitive.) 
  The homotopy category $\Ho(C)= C[W^{-1}]$ of a coloured category
  $(C,W)$ is the category obtained by formally inverting the equiarrows
  (cf.~\cite{Gabriel-Zisman}).
  Let $\CCat$ denote the category of
  coloured categories and colour-preserving functors (i.e.~functors
  preserving coloured arrows).
  
  Key examples of coloured categories are $\Top$ and $\sSet$,
  colouring the weak homotopy equivalences, and also $\Cat$ and
  $\nCat$, with the appropriate notions of fully faithful and
  essentially surjective functors as equimorphisms.  We shall also
  have good use of the coloured category $\Set$ in which the
  equimorphisms are the bijections.  (Note that any category admits
  three `trivial' colourings: colouring the identity arrows only,
  colouring the isomorphisms only, or colouring all arrows.)
\end{blanko}
  
\begin{blanko}{Tamsamani $\SS$-categories.}
%   We will work with a fixed coloured category $\SS$, and define 
%   categories in it.  We impose the extra axiom that $\SS$ has 
%   arbitrary sums and finite products, and we require the notion
%   of equimorphism to be stable under these operations.  In this 
%   situation, an object is called {\em contractible} if its morphism to the 
%   terminal object is an equimorphism.
% 
  Given a coloured category $\SS$ with discrete objects and admitting
  fibre products over discrete objects, it makes sense to define a
  {\em Tamsamani $\SS$-category} to be a functor $X:\Delta\op\to\SS$
  subject to the weak Segal condition and with $X_0$ discrete.  An
  important example of this generalisation is given by taking
  $\SS=\sSet$: these are called {\em Segal categories} 
  (cf.~\cite{Dwyer-Kan-Smith:1989}), and in many
  context they play the r\^ole of certain $\infty$-categories,
  cf.~Hirschowitz-Simpson~\cite{Hirschowitz-Simpson},
  To\"en-Vezzosi~\cite{Toen-Vezzosi:0212}, and
  To\"en~\cite{Toen:0409}.
  \nocite{Toen-Vezzosi:0207} %\nocite{Toen-Vezzosi:0404}
\end{blanko}

\begin{blanko}{Getting the colours right, and discrete objects.}
  Returning to the inductive definition of Tamsamani $n$-category
  (\ref{Tams}), in order to make sense of the induction step, we must
  define what an equimorphism of $\SS$-categories is, and the subtle
  point is to get this definition right.  The subtlety can be observed
  already in the definition of equimorphism in $\Cat$.  The notion of
  equimorphism is weaker than the pointwise one (because an
  equimorphism of categories is not necessarily a bijection on
  objects), and it is stronger than the notion induced from weak
  homotopy equivalence in $\sSet$ (e.g., as categories $\ord 0$ and
  $\ord 1$ are not equivalent, but their nerves are weakly
  equivalent).  The condition `admitting a quasi-inverse' happens to
  be the correct notion for ordinary categories, but already for
  $2$-categories it is too strong.  The good general notion of
  equimorphism turns out to be `essentially surjective and fully
  faithful', provided these notions are well chosen.  Tamsamani's
  definition --- which is good, in view of his main theorem, quoted in
  \ref{Tams-thm} below --- relies on certain
  properties of discrete objects.  The axiomatisation of these
  properties must take into account that they should reproduce
  themselves under the induction step; the details are given in 
  the Appendix.
\end{blanko}

%%%%%%%%%%%%%%%%%%%%%%%%%%%%%%%%%%%%%%%%%%%%%%%%%%
\subsection*{Coloured ordinals}
%%%%%%%%%%%%%%%%%%%%%%%%%%%%%%%%%%%%%%%%%%%%%%%%%%

\begin{blanko}{Coloured ordinals.}
  Since nonempty finite ordinals are fundamental categories, it
  seems worthwhile to look at the corresponding coloured notion
  (cf.~\cite{Kock:grosdelta}).  Recall that ordinals can be seen as
  free categories on linearly ordered graphs (strings of arrows).  
  Consider the coloured version of these three notions:
  
  A {\em coloured graph} is a graph some
  of whose edges have been singled out as coloured.  The {\em free
  coloured category} on a coloured graph is defined by taking $C$ to 
  be the free
  category on the whole graph and taking $W$ to be the free
  category on the coloured part of the graph (including all vertices). 
  This means that in a free coloured category, the composite of two
  arrows is an equimorphism if and only if both arrows are
  equimorphisms.  
%   (In particular all invertible arrows are equimorphisms.)
  Finally a {\em (finite) coloured ordinal} is the free coloured 
  category on a (finite) linearly ordered coloured graph.

  So what it all boils down to is to take finite strings of arrows, some of
  which are coloured.  Let $\T$ denote the full subcategory of $\CCat$
  consisting of the (finite and non-empty) coloured ordinals.  There are
  many other descriptions of this category
  (cf.~\cite{Kock:grosdelta}): as the category of 
  epimorphisms in $\Delta$; as the category of planar trees
  of height $2$; as opposite to the category of subdivided finite strict
  intervals; as a Grothendieck construction or moduli space for all
  coherent colourings of finite ordinals, etc.

  The category $\T$ plays the same r\^ole for coloured categories
  as $\Delta$ does for plain categories.  For example, a coloured category $C$
  can be described in terms of its {\em coloured nerve},
\begin{eqnarray*}
  \T\op & \longrightarrow & \Set  \\
  K & \longmapsto & \Hom_{\CCat}(K,C) .
\end{eqnarray*}

\end{blanko}

\begin{blanko}{Graphical interpretation.}
   We represent the objects of a coloured ordinal $K$ as dots
  arranged in a column (just like the drawing for ordinals on page 
  \ref{delta-graphic}).  The ordinary
  arrows are not drawn (the order expressed by the column indicates
  everything).  The equimorphisms are drawn as a link (the direction of the
  arrow being expressed by the order in the column: the arrows go upwards). 

  The graphical expression for functors between coloured ordinals is just
  as for usual ordinals what the dots are concerned, and for the links the
  rule is that a link can be set but may not be broken.  Here is a list of the
  most basic arrows in $\T$ (not mentioning the identity arrows):
\begin{center}\label{11}
\begin{texdraw}
	\setunitscale 0.75
	\move (-10 0)
	% THIS IS IDENTITY
	\bsegment
	\move (0 -5) \Bdot\hop\Bdot
	\move (20 0) \Bdot
	\linewd 0.2 \move (0 -5) \lvec (20 0) \lvec (0 5)
	\esegment
	\move (40 0)
	% THIS IS U-IDENTITY
	\bsegment
	\move (0 -5) \Bdot\Blink\Bdot
	\move (20 0) \Bdot
	\linewd 0.2 \move (0 -5) \lvec (20 0) \lvec (0 5)
	\esegment
	\move (120 0)
	% THIS IS U -> A
	\bsegment
	\move (0 -5) \Bdot\hop\Bdot
	\move (20 -5) \Bdot\Blink\Bdot
	\linewd 0.2 \move (0 -5) \lvec (20 -5) 
		    \move (0 5) \lvec (20 5)
	\esegment
	\move (200 0)
	% THIS IS SOURCE
	\bsegment
	\move (0 0) \Bdot
	\move (20 -5) \Bdot\hop\Bdot
	\linewd 0.2 \move (0 0) \lvec (20 -5)
	\esegment
	\move (250 0)
	% THIS IS TARGET
	\bsegment
	\move (0 0) \Bdot
	\move (20 -5) \Bdot\hop\Bdot
	\linewd 0.2 \move (0 0) \lvec (20 5)
	\esegment
	\move (300 0)
	% THIS IS SOURCE
	\bsegment
	\move (0 0) \Bdot
	\move (20 -5) \Bdot\Blink\Bdot
	\linewd 0.2 \move (0 0) \lvec (20 -5)
	\esegment
	\move (350 0)
	% THIS IS TARGET
	\bsegment
	\move (0 0) \Bdot
	\move (20 -5) \Bdot\Blink\Bdot
	\linewd 0.2 \move (0 0) \lvec (20 5)
	\esegment
	\move (420 -5)
	% THIS IS COMPOSITION
	\bsegment
	\move (0 0) \Bdot\hop\Bdot
	\move (20 -5) \Bdot\hop\Bdot\hop\Bdot
	\linewd 0.2 \move (0 0) \lvec (20 -5)
		    \move (0 10) \lvec (20 15)
	\esegment
	\move (470 -5)
	% THIS IS UA-COMPOSITION
	\bsegment
	\move (0 0) \Bdot\hop\Bdot
	\move (20 -5) \Bdot\Blink\Bdot\hop\Bdot
	\linewd 0.2 \move (0 0) \lvec (20 -5)
		    \move (0 10) \lvec (20 15)
	\esegment
	\move (520 -5)
	% THIS IS AU-COMPOSITION
	\bsegment
	\move (0 0) \Bdot\hop\Bdot
	\move (20 -5) \Bdot\hop\Bdot\Blink\Bdot
	\linewd 0.2 \move (0 0) \lvec (20 -5)
		    \move (0 10) \lvec (20 15)
	\esegment
	\move (570 -5)
	% THIS IS UU-COMPOSITION
	\bsegment
	\move (0 0) \Bdot\Blink\Bdot
	\move (20 -5) \Bdot\Blink\Bdot\Blink\Bdot
	\linewd 0.2 \move (0 0) \lvec (20 -5)
		    \move (0 10) \lvec (20 15)
	\esegment
\end{texdraw}
\end{center}
If $C:\T\op\to\Set$ is (the nerve of) a coloured category $(C,W)$, then the 
images of these basic arrows have precise interpretation as structure
maps, just as remarked in \ref{delta-graphic}: the first one 
associates the identity arrow to an object; ditto for the second but 
it specifies  that this arrow is an equimorphism; the third is the
inclusion of $W$ into $C$, and so on.  The second-to-last is the 
composition of an arbitrary arrow with an equimorphism.
\end{blanko}

\begin{blanko}{The projection $\T\to\Delta$ and colour structure on 
  $\T$.}
  Consider the natural projection functor $\pi:\T\to\Delta$ given by
  `taking equi-connected components', i.e., contracting all links.  Now 
  $\T$ has a natural colour structure, given by taking the
  equimorphisms $\T$ to be those arrows mapping to identity
  arrows in $\Delta$.  
%   This notion coincides with various different 
%   notions of equimorphism inherited from $\CCat$.
\end{blanko}
\nocite{Dwyer-Kan:simplicial-localization}
% \pagebreak

%%%%%%%%%%%%%%%%%%%%%%%%%%%%%%%%%%%%%%%%%%%%%%%%%%
\section{Semi-categories, coloured semi-categories, and the fat delta}
%%%%%%%%%%%%%%%%%%%%%%%%%%%%%%%%%%%%%%%%%%%%%%%%%%

We want to weaken the identity arrow axiom, but there is no way to do
that within the purely simplicial viewpoint: the identity arrows arise
inherently as a consequence of the degeneracy maps in $\Delta$. 
Plainly removing degeneracy maps leads to
$\Dmono$-diagrams, 
as we proceed to explain.  $\Dmono$-diagrams 
satisfying the Segal condition 
are semi-categories --- note that the Segal condition relates only to 
face maps.
Passing from categories to semi-categories is too drastic a reduction however,
and there are many constructions with categories that fail for
semi-categories.  The fat delta $\gros$ will be a sort of intermediate
between $\Dmono$ and $\Delta$.

\begin{blanko}{Semi-categories and `semi-ordinals'.}
  We use the prefix {\em semi-} consistently to mean `non-unital': A
  {\em semi-category} is just like a category, except that identity
  arrows are not required.  A {\em semi-functor} is a map compatible
  with the composition law.  Note that the identity semi-functor
  exists for any semi-category, so semi-categories and semi-functors
  form a genuine category $\semiCat$, not just a semi-category.  A
  (finite) {\em `semi-ordinal'} is the semi-category associated to a
  (finite) total strict order relation $<$.  Since $<$ is not
  reflexive, a given element is not related to itself, so there are no
  identity arrows.  As a consequence, all morphisms between
  semi-ordinals are injective, and the category of finite non-empty
  semi-ordinals is naturally identified with $\Dmono$.  The (semi-)nerve of a
  semi-category is a functor $\Dmono\op\to\Set$ satisfying the
  (strict) Segal condition.

  It is straightforward to copy over the definition of higher
  categories to the case of semi-categories. 
  But while the basic definitions of semi-categories are easy, their
  theory is quite different from that of categories.  For instance, in
  the theory of categories constant diagrams always exist, and one can
  define limits and colimits as adjoints of the constant diagram
  functor~\cite{MacLane:categories}.  For semi-categories, constant
  diagrams do not in general exist, because the supporting object is
  not required to have an identity arrow.
\end{blanko}

\begin{blanko}{Coloured semi-ordinals.}
  Combining all the previous notions we finally come to the promised
  fat delta: the category $\gros$ of coloured finite non-empty
  semi-ordinals.  The definitions should be obvious: a (finite) {\em
  coloured semi-category} is a semi-category with a sub-semi-category
  comprising all objects, and a morphism between coloured 
  semi-categories is a semi-functor required to preserve
  colour.  There is a (genuine) category $\kat{C\textonehalf Cat}$ of
  coloured semi-categories and their morphisms (and in fact this
  category can be considered a coloured category).
 
  The {\em coloured (finite) semi-ordinals} are the free coloured
  semi-categories on (finite) linearly ordered coloured graphs.  So it
  boils down to giving a (finite) string of arrows, some of which are
  coloured;  there are no identity arrows at the vertices.  
\end{blanko}
  
\begin{blanko}{The fat delta $\gros$.} 
  The fat delta $\gros$ is by definition the full subcategory of
  $\kat{C\textonehalf Cat}$ consisting of all finite non-empty
  coloured semi-ordinals.  It is naturally identified with the
  category of monomorphisms in $\T$:
  $$
  \gros = \T_{\operatorname{mono}} .
  $$
  The various characterisations of $\T$ yield alternative descriptions of
  $\gros$.  In particular we get the following concise description of
  $\gros$ as a subcategory of the
  category of arrows in 
$\Delta$: its objects are the epimorphisms in $\Delta$, and the
arrows are the monomorphisms in $\Arr(\Delta)$;  these are the 
commutative diagrams in $\Delta$ whose downward arrows are 
epimorphisms and whose top arrow is a monomorphism:
\begin{diagram}[w=3.6ex,h=3.6ex,tight]
\cdot & \rInto  & \cdot  \\
\dEpi  &    & \dEpi  \\
\cdot  & \rTo  & \cdot
\end{diagram}

  The drawings of objects and arrows in $\gros$ are the same as for 
  $\T$, except that the maps are required to
  be `injective on dots'.  Thus the first two figures listed on
  page~\pageref{11} are not arrows in $\gros$.
\end{blanko}
  
\begin{blanko}{The projection $\gros\to\Delta$ and colours in $\gros$.}
  An important r\^ole is played by the projection functor
  $\pi:\gros\to\Delta$, sending a coloured (semi-) ordinal to the
  ordinal of equi-connected components; let $V\subset \gros$ denote the
  subcategory of vertical arrows for $\pi$ (i.e.~whose image in 
  $\Delta$ is an identity arrow), then the pair $(\gros,V)$
  is a coloured category.  
%   This notion of colour coincides with
%   various notions induced from $\kat{C\textonehalf Cat}$.  
  From the drawing on page~\pageref{11}, these are equimorphisms:
  \begin{center}
  \begin{texdraw}
	  \setunitscale 1
	  \move (0 0)
	  % THIS IS SOURCE
	  \bsegment
	  \move (0 0) \Bdot
	  \move (20 -5) \Bdot\Blink\Bdot
	  \linewd 0.2 \move (0 0) \lvec (20 -5)
	  \esegment
	  \move (50 0)
	  % THIS IS TARGET
	  \bsegment
	  \move (0 0) \Bdot
	  \move (20 -5) \Bdot\Blink\Bdot
	  \linewd 0.2 \move (0 0) \lvec (20 5)
	  \esegment
	  \move (130 -5)
	  % THIS IS UA-COMPOSITION
	  \bsegment
	  \move (0 0) \Bdot\hop\Bdot
	  \move (20 -5) \Bdot\Blink\Bdot\hop\Bdot
	  \linewd 0.2 \move (0 0) \lvec (20 -5)
		      \move (0 10) \lvec (20 15)
	  \esegment
	  \move (180 -5)
	  % THIS IS AU-COMPOSITION
	  \bsegment
	  \move (0 0) \Bdot\hop\Bdot
	  \move (20 -5) \Bdot\hop\Bdot\Blink\Bdot
	  \linewd 0.2 \move (0 0) \lvec (20 -5)
		      \move (0 10) \lvec (20 15)
	  \esegment
	  \move (260 -5)
	  % THIS IS UU-COMPOSITION
	  \bsegment
	  \move (0 0) \Bdot\Blink\Bdot
	  \move (20 -5) \Bdot\Blink\Bdot\Blink\Bdot
	  \linewd 0.2 \move (0 0) \lvec (20 -5)
		      \move (0 10) \lvec (20 15)
	  \esegment
  \end{texdraw}
  \end{center}
and these are not:
\begin{center}
\begin{texdraw}
	\setunitscale 1
	\move (0 0)
	% THIS IS U -> A
	\bsegment
	\move (0 -5) \Bdot\hop\Bdot
	\move (20 -5) \Bdot\Blink\Bdot
	\linewd 0.2 \move (0 -5) \lvec (20 -5) 
		    \move (0 5) \lvec (20 5)
	\esegment
	\move (80 0)
	% THIS IS SOURCE
	\bsegment
	\move (0 0) \Bdot
	\move (20 -5) \Bdot\hop\Bdot
	\linewd 0.2 \move (0 0) \lvec (20 -5)
	\esegment
	\move (130 0)
	% THIS IS TARGET
	\bsegment
	\move (0 0) \Bdot
	\move (20 -5) \Bdot\hop\Bdot
	\linewd 0.2 \move (0 0) \lvec (20 5)
	\esegment
	\move (210 -5)
	% THIS IS COMPOSITION
	\bsegment
	\move (0 0) \Bdot\hop\Bdot
	\move (20 -5) \Bdot\hop\Bdot\hop\Bdot
	\linewd 0.2 \move (0 0) \lvec (20 -5)
		    \move (0 10) \lvec (20 15)
	\esegment
\end{texdraw}
\end{center}

These last four arrows should be compared to the four arrows in $\Delta$ 
drawn on page~\pageref{4}.  
The important thing to note, compared to the situation in $\Delta$ is 
that the figure 
\raisebox{-2pt}{
\begin{texdraw}	\setunitscale 0.75
  % THIS IS IDENTITY
  \move (-3 0)
  \bsegment
  \move (0 -5) \Bdot\hop\Bdot
  \move (20 0) \Bdot
  \linewd 0.2 \move (0 -5) \lvec (20 0) \lvec (0 5)
  \esegment
  \move (23 0)
\end{texdraw}
}
(which was responsible for the existence of identity arrows), does 
{\em not} exist in $\gros$.  Instead we have 
\raisebox{-2pt}{
\begin{texdraw} \setunitscale 0.75
  % THIS IS IDENTITY
  \move (-3 0)
  \bsegment
  \move (0 -5) \Bdot\hop\Bdot
  \move (20 -5) \Bdot\Blink\Bdot
  \linewd 0.2 \move (0 -5) \lvec (20 -5)\move (0 5) \lvec (20 5)
  \esegment
  \move (23 0)
\end{texdraw}}.
\end{blanko}

Here is a picture of the first few arrows in $\grosdelta$:
\begin{diagram}[w=5ex,h=5ex,tight,pilespacing=5pt]
\intextB & \pile{\rTo \\ \rTo}  & \intextBB  \\
\dTo\dTo  & \ldTo   &   \\
\intextBbB  &   & 
\end{diagram}
and their images in $\Delta$:
\begin{diagram}[w=5ex,h=5ex,tight,pilespacing=5pt]
\intextB & \pile{\rTo \\ \rTo \\ \lTo[shortfall=2ex]}  & \intextBB 
\end{diagram}

\begin{blanko}{The two inclusions $\Dmono\subset \gros$.}
  In the triangle diagram above we see the beginning of two copies
  of $\Dmono$.  There is the `horizontal' inclusion
  $\iota:\Dmono\to \gros$, interpreting a semi-ordinal as a coloured 
  semi-ordinal with nothing coloured:
  \begin{diagram}[w=5ex,h=5ex,tight,pilespacing=5pt]
  \intextB & \pile{\rTo \\ \rTo}  & \intextBB   & \pile{\rTo \\ \rTo 
  \\ \rTo} & \intextBBB & \ldots
  \end{diagram}
  The composite functor
  $$
  \Dmono \into \gros \twoheadrightarrow \Delta
  $$
  is just the standard inclusion of the monos in $\Delta$.  This is
  the sense in which the fat delta is intermediate  between $\Dmono$ and
  $\Delta$.
  
  There is also the `vertical' inclusion $\zeta:\Dmono \into \gros$, 
  interpreting a semi-ordinal as a coloured 
  semi-ordinal with everything coloured.  The image is
  \begin{diagram}[w=5ex,h=5ex,tight,pilespacing=5pt]
  \intextB & \pile{\rTo \\ \rTo}  & \intextBbB   & \pile{\rTo \\ \rTo 
  \\ \rTo} & \intextBbBbB & \ldots
  \end{diagram}
  (for typographical reasons drawn horizontally, but to fit into the 
  drawing above it should be vertical).
\end{blanko}

\begin{blanko}{Discrete object structure on $\gros$ (and $\grosop$).}
  Just like for $\Delta$ (and $\Delta\op$), we declare the single dot
  to be the only discrete object in $\gros$ (or in $\grosop$).
%   Thinking of `discreteness' as dual to `connectedness', this is the
%   only meaningful definition of discrete objects in $\gros$, although
%   naturally it is not covered by the notion of standard discrete
%   objects introduced in the Appendix (\ref{SDO}),
%   designed for categories with better closure properties.
\end{blanko}

% \pagebreak
%%%%%%%%%%%%%%%%%%%%%%%%%%%%%%%%%%%%%%%%%%%%%%%%%%
\section{Definition of fair categories and fair $n$-categories}
%%%%%%%%%%%%%%%%%%%%%%%%%%%%%%%%%%%%%%%%%%%%%%%%%%

%%%%%%%%%%%%%%%%%%%%%%%%%%%%%%%%%%%%%%%%%%%%%%%%%%
\subsection*{Fair $\SS$-categories}
% %%%%%%%%%%%%%%%%%%%%%%%%%%%%%%%%%%%%%%%%%%%%%%%%%%

Henceforth, let $\SS$ denote a coloured category with a notion of
discrete objects.  (It is reasonable to assume the existence of fibre
products over discrete objects, and require this operation to preserve
equiarrows, but in fact the definition makes sense without these
requirements.)

\begin{blanko}{The idea.}\label{idea}
  A fair category in $\SS$ is going to be a 
  $\grosop$-diagram
  $$
  X : \grosop \to \SS
  $$
  satisfying three obvious axioms.  In order to exhibit the axioms as
  obvious, let us right away explain how $X$ is to be interpreted as a 
  nerve:
  put
  $$
  O\df X\!_{\intextB}
  \qquad
  A \df X\!_{\intextBB}
  \qquad
  U \df X\!_{\intextBbB} ,
  $$
  and think of these as the spaces 
  of {\em objects}, {\em arrows}, and {\em (weak) identity arrows}, respectively.
  This notation will be used throughout.
  % We will refer to an $\SS$-category by the triple $(O,A,U)$.
  
  The images of the first few maps (i.e., the beginning of the $\grosop$-diagram)
  looks like this:
  \begin{equation}\label{triangle}
  \begin{diagram}[w=5ex,h=5ex,tight,pilespacing=5pt]
    O & \pile{\lTo^t \\ \lTo_s}  & A  \\
    \uTo\uTo  & \ruTo_u    &   \\
    U  &   & 
  \end{diagram}
  \end{equation}
  The maps $s,t$ are {\em source} and {\em target}, and $u$ is thought of as
  the inclusion of the space of identity arrows into the space of all arrows.
\end{blanko}

\begin{blanko}{Discreteness.}
  Now the first axiom is clear:
  \begin{quote}
    The object $O$ (image of the single dot \intextB) is discrete.
  \end{quote}
  Taking \intextB\ as the only discrete object in $\grosop$, the first
  axiom says that $X: \grosop\to\SS$ preserves discrete objects.
\end{blanko}

\begin{blanko}{Segal condition.}
  Second, we want a Segal condition which should induce the usual 
  Segal condition on each  $\Dmono\subset\gros$.  Just like in $\Delta$ 
  (and in $\Dmono$), pushouts over the single dot exists in $\gros$
  if and only if
  the dot is included as last dot in the first summand and as the 
  first dot in the second (or vice versa).
  Let the pushout of $m \lTo \intextB \rTo
  n$ be denoted by $m\dotsum n$.
  The {\em generalised Segal condition} requires the
  {\em generalised Segal maps}
  $$
  X_{ m \smalldotsum n } \rTo X_m \times_O X_n
  $$
  to be isomorphisms.  This formulation presupposes
  existence of fibre products over discrete objects, but
  we might as well just say that for all $m$ and $n$, the square
  \begin{diagram}[w=6ex,h=4.5ex,tight]
  X_{ m \smalldotsum n } & \rTo  & X_n  \\
  \dTo  &    & \dTo  \\
  X_m  & \rTo  & O
  \end{diagram}
  should be a pullback square. 
  Hence our second axiom is: 
  \begin{quote}
    $X$ preserves fibre products over discrete objects.
  \end{quote}

  In practice, this means two things.  First of all, the restriction
  to either copy of $\Dmono\subset \gros$ is a $\Dmono\op$-diagram
  which satisfies the Segal condition.  Hence $A$ and $U$ are each
  semi-categories, i.e.~carry associative composition operations over
  $O$.  The $\Dmono\op$-diagram
  $\Dmono\op\stackrel{\iota}{\longrightarrow} \grosop
  \stackrel{X}{\longrightarrow} \SS$ is the {\em underlying
  semi-category of $X$.} 
  
  Second, the Segal condition means that that
  the rest of the diagram can be constructed from the above triangle
  diagram (\ref{triangle}), provided $U \to A$ is a semi-functor 
  (i.e.~is compatible with composition in $U$ and in $A$).
\end{blanko}

\begin{blanko}{Weak identity arrows.}\label{colour-axiom}
  Finally we want the points in $U$ to act as identity arrows.  For each
  object there should be a contractible space of weak identity arrows, so
  we want equimorphisms for the two maps $U \topile O$, images of
   \begin{center}
  \begin{texdraw}
	  \setunitscale 1
	  \move (0 0)
	  % THIS IS SOURCE
	  \bsegment
	  \move (0 0) \Bdot
	  \move (20 -5) \Bdot\Blink\Bdot
	  \linewd 0.2 \move (0 0) \lvec (20 -5)
	  \esegment
	  \move (50 0)
	  % THIS IS TARGET
	  \bsegment
	  \move (0 0) \Bdot
	  \move (20 -5) \Bdot\Blink\Bdot
	  \linewd 0.2 \move (0 0) \lvec (20 5)
	  \esegment
  \end{texdraw}
  \end{center}
  The identity condition is that composing with a weak identity arrow
  should be neutral up to homotopy.  So we want the composition maps
  $U \times_O A \rTo A \lTo A \times_O U$  and $U \times_O U \rTo U$
  to be 
  equimorphisms --- these are the image of 
  \begin{center}
  \begin{texdraw}
	  \setunitscale 1
	  \move (130 -5)
	  % THIS IS UA-COMPOSITION
	  \bsegment
	  \move (0 0) \Bdot\hop\Bdot
	  \move (20 -5) \Bdot\Blink\Bdot\hop\Bdot
	  \linewd 0.2 \move (0 0) \lvec (20 -5)
		      \move (0 10) \lvec (20 15)
	  \esegment
	  \move (180 -5)
	  % THIS IS AU-COMPOSITION
	  \bsegment
	  \move (0 0) \Bdot\hop\Bdot
	  \move (20 -5) \Bdot\hop\Bdot\Blink\Bdot
	  \linewd 0.2 \move (0 0) \lvec (20 -5)
		      \move (0 10) \lvec (20 15)
	  \esegment
	  \move (260 -5)
	  % THIS IS UU-COMPOSITION
	  \bsegment
	  \move (0 0) \Bdot\Blink\Bdot
	  \move (20 -5) \Bdot\Blink\Bdot\Blink\Bdot
	  \linewd 0.2 \move (0 0) \lvec (20 -5)
		      \move (0 10) \lvec (20 15)
	  \esegment
  \end{texdraw}
  \end{center}
respectively.  
  Note that the five maps mentioned 
  here are all images of
  vertical maps in $\grosdelta$.  In fact, these five maps generate the
  category of vertical arrows in $\gros$ under the $\dotsum$ operation
  (see \cite{Kock:grosdelta} for details), so in view 
  of the Segal condition and the assumption that equimorphisms in 
  $\SS$ are stable under fibre products over discrete objects, the 
  requirement that these five maps be equimorphisms is equivalent to
  requiring {\em every}
  vertical map in $\grosdelta$ to be sent to an equimorphism in $\SS$.
  This is of course a more uniform and conceptual condition (and it 
  makes sense without the extra assumptions), so the third axiom is:
\begin{quote}
  $X:\grosop\to\SS$  preserves colours.
\end{quote}
\end{blanko}

In short, the definition is simply this:

\begin{blanko}{Definition of fair category.}\label{def}
  A {\em fair $\SS$-category} is a colour-preserving
  functor \mbox{$\grosop\to\SS$,} preserving also discrete objects and fibre
  products over discrete objects. 

  A {\em morphism} of fair $\SS$-categories is just
  a natural transformation.  Let $\SFairCat$ denote the category of fair
  $\SS$-categories and their morphisms.  Note that these morphisms are
  strict, and in some situations a weaker notion is needed.  For
  example, in order to get the correct inner homs, one
  should consider derived morphisms, just like in the theory of
  Tamsamani categories, cf.~\cite{Simpson:internalhom}.  See
  Remark~\ref{weak-morphisms} below for a particular case of a weaker
  notion of morphism.
\end{blanko}

\begin{blanko}{Strict categories, fair categories, and semi-categories.}
  \label{interpretations}
  There is an obvious forgetful functor $\SFairCat \to
  \kat{S-}\semiCat$, defined by restriction to $\Dmono$ via the
  horizontal inclusion $\iota: \Dmono\to\gros$.  Similarly every 
  strict $\SS$-category is canonically interpreted as a fair 
  $\SS$-category: pre-composition with $\pi: \gros\to\Delta$
  yields a full embedding $\pi\upperstar: \SCat \to \SFairCat$.
% Clearly this functor preserves discrete objects and 
% fibre products ovre discrete objects.
  Non-strict examples of fair $S$-categories are given in
  Section~\ref{Sec:examples}.
\end{blanko}

%%%%%%%%%%%%%%%%%%%%%%%%%%%%%%%%%%%%%%%%%%%%%%%%%%
\subsection*{Equimorphisms of fair $S$-categories}
%%%%%%%%%%%%%%%%%%%%%%%%%%%%%%%%%%%%%%%%%%%%%%%%%%

\begin{blanko}{Standard discrete objects.}
  Just as in the case of Tamsamani $\SS$-categories, in order to 
  get a good notion of equimorphism, further axioms must be imposed 
  on the notion of discrete objects and their relation to the 
  subcategory of equiarrows in $\SS$.  From now on,
  \begin{quote}\em
    We assume that
  $\SS$ has standard discrete objects with compatible colours,
  \end{quote}
  in the sense of the Appendix.
  In particular there is an adjunction
  $\pi_0 \isleftadjointto \delta$:
  $$
  \SS \pile{\lTo^\delta \\ \rTo_{\pi_0}} \Set  ;
  $$
  the objects in the image of $\delta$ are the discrete objects, and
  the left adjoint $\pi_0$ is the `components functor'.  Both these
  functors preserve discrete objects, and fibre products
  over discrete objects.  The compatibility of colours with respect to
  the discrete objects means that $\delta$ and $\pi_0$ are required to
  preserve equimorphisms.   Here and throughout, $\Set$ is considered a
  coloured category by taking the bijections as equimorphisms.
\end{blanko}

\begin{blanko}{Hom spaces.}\label{slice}
  It follows from the decomposition property (\ref{decomp})
  that in a fair $\SS$-category $X=(O,A,U)$ we have a
  decomposition of the space of arrows
  $$
  A = \coprod_{x,y\in O} A(x,y)
  $$
  according to source and target.   Here $A(x,y)$ is the subspace 
  of $A$ consisting of all arrows whose source
  is $x$ and whose target is $y$.
%   These are the hom spaces, in terms 
%   of which we shall define the notion of `fully faithful'.
\end{blanko}

\begin{lemma}\label{s=t}
  In a fair $\SS$-category $X=(O,A,U)$, the two maps $O \topileback U$
  coincide.  In other words, weak identity arrows are endomorphisms.
\end{lemma}

\begin{dem}
  The composite
  $$
  \grosop \stackrel{X}{\rTo} \SS \stackrel{\pi_0}{\rTo} \Set
  $$
  is a fair category in $\Set$, since $\pi_0$
  preserves discrete objects, fibre products over discrete objects,
  and equimorphisms (all the notions involved in the definition of
  fair category).  Since $O$ is discrete and $\pi_0$ is the 
  reflector, we have a bijection
  $$
  \SS(U,O)
%   \simeq \Set(\pi_0 U, O)
  \leftrightarrow \Set (\pi_0 U , \pi_0 O),
  $$
%   expressing that $\Set$ is a reflexive subcategory in $\SS$, 
  so it is enough to establish the result in the case of a fair category
  in $\Set$.  The two maps $O \topileback U$ are the beginning of a
  $\Dmono\op$-diagram
  $$
  \newdiagramgrid{faces}{1,1,1,1.5}{1}
  \begin{diagram}[w=6ex,h=4.5ex,grid=faces,tight]
  O & \pile{\lTo^s\\ \lTo_t}  
  & U &  \pile{\lTo^a\\ \lTo~b \\ \lTo_c} 
  & U \times_O U & \dots
  \end{diagram}
  $$
  and since we are in $\Set$, all these maps are bijections.  The
  face map identities read $s\circ a = s\circ b$, $t\circ b = t\circ 
  c$, and $s\circ c = 
  t\circ a$.  Since $s$ and $t$ are invertible, the first two equations
  imply $a=b=c$.  Now the third equation can be written $s\circ b = 
  t\circ b$, and since $b$ is invertible we conclude $s=t$. 
%   (In fact 
%   the same argument shows more generally that any two parallel arrows in a 
%   $\Dmono\op$-diagram in a groupoid coincide.)
\end{dem}

\bigskip

The definition of fair $n$-category will be inductive: a fair 
$n$-category is a fair category in the coloured category of 
fair $(n-1)$-categories.
The base for this induction, and the key to understanding
$\grosop$-diagrams as categories, is the case where $\SS$ is
the category of sets, with bijections as coloured arrows.
% Several categorical notions will be defined in terms of comparison
% with this case.

\begin{blanko}{Fair $\Set$-categories.}\label{1FCat}
  We consider $\Set$ a coloured category by taking the bijections as
  equi\-arrows.  Given a fair category $\grosop\to\Set$, denoted
  $X=(O,A,U)$ as in \ref{idea}, then $U \to O$ is a bijection, and 
  therefore $U \to A$ is an injection, as seen in the diagram in 
  \ref{idea}.  The restriction of $X$ to
  $\iota:\Dmono\op\into\grosop$ is a semi-category $O \topileback A$ with
  the property that every object has an identity arrow.  This follows
  because the inverse of $U \to O$, followed by $U \to A$, provides
  the degeneracy map $O \To A$ (and the remaining degeneracy maps are 
  provided by the Segal condition); the commutativity of the
  $\grosop$-diagram immediately implies that the degeneracy map
  identities are satisfied.  This construction defines a functor
  $$
  \theta : \kat{Set-FairCat} \rTo \Cat .
  $$
\end{blanko}

\begin{blanko}{Fair nerve of a category.}\label{fairnerve}
  In the other direction, starting from any category $\CC$, there is
  a natural `fair nerve' associated with it: let $\CC$ be coloured by
  taking the identity arrows as equiarrows, then the fair nerve 
  functor
  $$
  \rho : \Cat \to \kat{Set-FairCat}
  $$
  is defined by sending a category 
  $\CC$ to
 \begin{eqnarray*}
    \grosop & \longrightarrow & \Set  \\
    K & \longmapsto & \Hom_{\CCat}( K, \CC) ,
  \end{eqnarray*}
  which is readily seen to be a fair category in $\Set$.  If we denote it
  $X=(O,A,U)$ like in 
  \ref{idea}, then $O$ is the set of all objects in $\CC$, 
  $A$ is the set of all arrows in $\CC$, and
  $U$ is the set of all identity arrows in $\CC$. 
%   Let $i: U \to A$ be the inclusion of the identity arrows into the 
%   set of all arrows.
  In this interpretation, the image of the fair nerve functor consists
  of those fair $\Set$-categories $(O,A,U)$ for which $U \to A$ is an
  inclusion, not just an injection.  This is a full reflective
  subcategory; the reflector is $\theta$.

  Starting with a fair $\Set$-category $(O,A,U)$, applying $\theta$
  gives the semi-category $O \topileback A$ which happens to be a
  category, i.e.~has identity arrows $u(U)$ (and the morphisms happen
  to preserve these), and then applying $\rho$ gives the fair category
  $(O,A,u(U))$, so the unit for the adjunction is the isomorphism
  $(O,A,U) \to (O,A,u(U))$ given as the identity functor on the
  $A$-semi-nerve, and the isomorphism $u: U \to u(U)$ on the 
  $U$-semi-nerve.
  
  Conversely starting with a category $\CC = (O \topileback A)$,
  applying $\rho$ we first get the fair category $(O,A,U=\Id(\CC))$,
  and then we throw the $U$-part away and observe that the resulting
  semi-category is just the category $\CC$ again.

  In summary:
\end{blanko}
\begin{prop}\label{dim1}
  There is an adjoint equivalence of categories
$$
\kat{Set-FairCat} \pile{\lTo \\ \rTo} \Cat .
$$
\end{prop}

In the above back-and-forth construction we took a category 
viewpoint, focusing on the start of the $\grosop$-diagram $(O,A,U)$
and letting the Segal condition take care of the rest.  But in fact 
the same result holds without requiring the Segal condition:
\begin{prop}\label{sSet=FairSSet}
  There is an adjoint equivalence of categories
  $$
  \CCat(\grosop,\Set) \simeq \Cat(\Delta\op,\Set) = \sSet .
  $$
\end{prop}

For the fine points of the appendix, it is important to observe that
the functors of these two propositions preserve sums and finite
products, as well as discrete objects.

\bigskip

We shall give $\SFairCat$ the structure of a coloured
category, and then the above equivalence is an equivalence of coloured
categories.  Just as for ordinary categories, a functor $F :
A\rightarrow B$ of fair categories in $\SS$ will be called an
equimorphism if it is fully faithful and essentially surjective.
These notions must now be defined.

\begin{blanko}{Fully faithfulness.}\label{ff}
  Given a fair category $X=(O,A,U)$, write $A = \disju_{x,y\in O} A(x,y)$
  and $U = \disju_{x\in O} U(x)$, as in Remark~\ref{slice}.  Now a fair
  functor $F:X\Rightarrow X'$ is said to be {\em fully faithful} if for
  each pair of objects $x,y$ the map $A(x,y) \to A'(Fx, Fy)$ is an
  equimorphism in $\SS$, and for each object $x$ the map $U(x) \to U'(Fx)$ 
  is 
  an equimorphism.  (Note that this last condition is automatic if
  the two-out-of-three property holds for 
  equimorphisms in $\SS$.)
%   
%   Note  that this definition 
%   relies on the decomposition property \ref{decomp}.)
\end{blanko}

The notion of essential surjectivity is subtler, and depends on a
notion of truncation.  We take our clue from usual categories:
consider the truncation functor $\tau_0 : \Cat \To \Set$ which to a
category associates the set of isomorphism classes of its objects.
Note that this functor is different from the components functor
$\pi_0$, but that the two functors agree on the discrete categories.
A morphism in $\Cat$ is essentially surjective if its truncation is a
surjection of sets.  The crucial features of $\tau_0$ are that it
comes with a natural transformation $\tau_0 \Rightarrow \pi_0$, and
that it preserves finite products and equimorphisms.  We take these as
axioms for the truncation functor:

\begin{blanko}{Truncation.}\label{truncation}
  A {\em truncation functor} on a coloured category $\SS$ with
  standard discrete objects is a colour-preserving functor
  $\tau_0:\SS\to\Set$, equipped with a natural transformation $\tau_0
  \Rightarrow \pi_0$, required to preserve arbitrary sums and
  finite products.  Since sums are preserved, $\tau_0$ and
  $\pi_0$ agree on discrete objects.  It follows from the
  decomposition property (\ref{decomp})
  that $\tau_0$ also preserves fibre products over
  discrete objects.
% \begin{MOREDETAILS}
%   Idea of proof.  Given
%   \begin{diagram}[w=6ex,h=4.5ex,tight]
%   P  & \rTo  & Y  \\
%   \dTo  &    & \dTo  \\
%   X  & \rTo  & D
%   \end{diagram}
%   where $D=\coprod_{s\in S} *$, write (check this!) $P = \coprod_{s\in S} P(s)
%   = \coprod_{s\in S} (X(s) \times Y(s))$.  Since $\pi_0$ preserves 
%   products it also preserves this fibre product over the discrete 
%   object $D$
% \end{MOREDETAILS}
\end{blanko}

\begin{blanko}{Essential surjectivity.}
  For a fixed  truncation functor  $\tau_0:\SS\to\Set$, there is 
  induced a functor 
  \begin{equation}\label{taustar}
  \SFairCat \rTo^{\tau_0{}\lowerstar} \Set\kat{-FairCat} \simeq \Cat .
  \end{equation}
  by sending a $\grosop\To\SS$ to $\grosop\To\SS\rTo^{\tau_0}\Set$ 
  and invoking Proposition~\ref{dim1}.  This works because $\tau_0$ 
  preserves the notions involved in the definition of fair category.
  
  Now a morphism $F: X \to X'$ in $\SFairCat$ is said to be {\em
  essentially surjective} (relative to $\tau_0$)
  if $\tau_0{}\lowerstar F$ is essentially surjective in
  the ordinary sense (it is a functor between ordinary categories).
  It should be noted that this notion of essentially surjective might
  not be the correct one, unless combined with fully faithful.  But
  that situation is all we care about:
%   Equivalently, it is essentially surjective if the long composite
%   is a surjective mapping of sets.
\end{blanko}

\begin{blanko}{Equivalences of fair categories.}\label{equiv}
  Let $\SS$ be a coloured category with standard discrete objects
  and a truncation functor $\tau_0:\SS\to\Set$.
  A morphism of fair $\SS$-categories is called an {\em equimorphism} if
  it is fully faithful and essentially surjective (with respect to 
  $\tau_0$).  Henceforth this is the notion referred to when talking
  about the coloured category $\SFairCat$.  Note that equimorphisms
  between 
  strict $\SS$-categories (considered as fair $\SS$-categories 
  via $\pi\upperstar : \SCat \to \SFairCat$, 
  cf.~\ref{interpretations}) are precisely the usual equivalences
  of $\SS$-categories.
\end{blanko}

\begin{lemma}
  Equimorphisms in $\SFairCat$ are stable under sums and finite 
  products.
\end{lemma}

\begin{dem}
  In fact this is true independently for 'fully faithful' and 'essential 
  surjective'.  For 'fully faithful' it follows because the hom `sets'
  of a sum (resp.~a finite product) is the sum (resp.~the product) of the hom
  `sets', and equimorphisms in $\SS$ are stable under sums 
  (resp.~finite products).  For `essentially surjective' it follows 
  because $\tau_0$ preserves sums and finite products.  
\end{dem}

\begin{prop}
  $\SFairCat$ has standard discrete objects.
\end{prop}
This is Proposition~\ref{SFairCat-SDO}.  The discrete-objects 
adjunction is given by
  $$
  \SFairCat \pile{\lTo^{\delta\lowerstar} \\ \rTo_{\pi_0{}\lowerstar }}
  \Set\kat{-FairCat} \simeq \Cat \pile{\lTo \\ \rTo}  
  \Set ,
  $$
where $\pi_0{}\lowerstar $ is defined by sending a fair $\SS$-category
$\grosop\To\SS$ to $\grosop\To\SS\rTo^{\pi_0}\Set$.

\begin{lemma}
  The components functor $\SFairCat \to \Set$ preserves equimorphisms
  as defined in \ref{equiv}.
\end{lemma}
  
\begin{dem}
  Note first that since $\pi_0: \SS\to\Set$ preserves equimorphisms,
  and since the definitions of fully faithful in $\SFairCat$ and
  $\kat{Set-FairCat}$ are both defined in terms of hom sets (in $\SS$
  in $\Set$ respectively), it follows that $\pi_0{}\lowerstar $
  preserve fully faithfulness.
  
%   If $F: X\Rightarrow X'$ is an essentially surjective functor of
%   categories $X,X' : \Delta\op\to \Set$, then $\delta\lowerstar F$ is
%   essentially surjective too, because its truncation
%   $\tau_0(\delta\lowerstar F)$ is just $F$ composed with $\tau_0 \circ
%   \delta$ which is isomorphic to $F$ itself since $\tau_0$ is a
%   retraction of $\delta$.  So $\delta\lowerstar $ preserves
%   equimorphisms.
   
  Seeing that $\pi_0{}\lowerstar$ preserves essential surjectivity
  relies on the natural transformation $u:\tau_0
  \Rightarrow \pi_0$.  
 For a fair $\SS$-category $X : \grosop \to \SS$, the 
 two categories $\tau_0{}\lowerstar  X$ and $\pi_0{}\lowerstar 
  X$ have the same object set, since both $\tau_0$ and $\pi_0$
  are the identity on discrete objects. 
  To say that a morphism $F: X \to X'$ of fair
  $\SS$-categories is essentially surjective means that 
  $\tau_0{}\lowerstar F$
  is an essentially surjective functor of categories.  This in turn 
  means that
  for every object $x'$ of $\tau_0{}\lowerstar X'$ 
  there exists an object $x$ of $\tau_0{}\lowerstar X$ and an
  isomorphism $\phi\in \tau_0{}\lowerstar  X'(Fx,x')$.
  But then $u(\phi)$ is an
  isomorphism in $\pi_0{}\lowerstar  X'(Fx, x')$ witnessing that 
  $\pi_0{}\lowerstar
  F$ is essentially surjective too, as required.
\end{dem}
  
Combining the three previous results we find that:
\begin{prop}\label{equiv-ditto}
  $\SFairCat$ has standard discrete objects
  with compatible colours.
\end{prop}

\begin{blanko}{Truncation for $\SFairCat$.}
  Given a truncation functor $\tau_0 : \SS\to\Set$, define a
  truncation functor $\SFairCat \to \Set$ as the composite
  $$
  \SFairCat \rTo^{\tau_0{}\lowerstar} \Set\kat{-FairCat} 
  \simeq \Cat \rTo^{\tau_0^{\Cat}} \Set .
  $$
  It follows readily that this new truncation functor preserves sum 
  and finite products (and it preserves equimorphisms by construction).
  The natural transformation to the components functor is the 
  horizontal composition of $u\lowerstar : \tau_0{}\lowerstar 
  \Rightarrow \pi_0{}\lowerstar $ and the natural transformation
  $\tau_0^{\Cat} \Rightarrow \pi_0^{\Cat}$.
\end{blanko}

Now $\SFairCat$ has been equipped with the same type of structure as $\SS$,
and the induction works:

\begin{blanko}{Fair $n$-categories.}\label{fair-n}
  A {\em fair $0$-category} is defined to be just a set.  
  Assuming we have
  already defined the coloured category $\kat{(n-1)FairCat}$ of all fair
  $(n-1)$-categories, we define a {\em fair $n$-category} to be a
  colour-preserving functor
  $$
  \grosop\to \kat{(n-1)FairCat}
  $$
  preserving discrete objects and fibre products over discrete objects
  (this is the strict Segal condition).
  Let $\kat{nFairCat}$ be the coloured category whose objects are the
  fair $n$-categories, whose morphisms are the natural 
  transformations,
  and the equimorphisms are as defined in \ref{equiv}.
  
  Notice that a strict $n$-category can be regarded as a fair
  $n$-category in a canonical way.  This follows by induction from
  \ref{interpretations} and the observation that the embedding
  $\kat{(n-1)Cat} \to \kat{(n-1)FairCat}$ preserves discrete objects, fibre 
  products over discrete objects, and equimorphisms (cf.~\ref{equiv}).
\end{blanko}

\begin{blanko}{Fair Tamsamani $n$-categories.}
  Repeating the definition of fair category (and the other definitions
  involved), but with weak Segal condition instead of the strict one,
  yields a notion of {\em fair Tamsamani $n$-categories}.  An
  interesting example of a fair Tamsamani $2$-category is given on
  page~\pageref{sub:Cob}.  Fair Tamsamani $n$-categories should be the
  context for comparing Tamsamani $n$-categories with fair
  $n$-categories, since both can be regarded as special cases of this
  notion.  It is expected that given such a fair Tamsamani
  $n$-category then there exists an equivalent Tamsamani
  $n$-category (strict identities).  In other words, weak identities
  can be strictified in the context of fair Tamsamani categories.
  Also there should exist an equivalent fair $n$-category, i.e., weak
  composition laws can be strictified in the context of fair Tamsamani
  categories. This last statement is one
  version of Simpson's conjecture~\ref{Simpson-conj}.
%   
%   Key ingredients in a proof of the second  result  would be 
%   Conjecture~\ref{Lgros} and a result of Dwyer-Kan~\cite{Dwyer-Kan:Princeton}
%   on homotopy types of restricted diagrams.
  
% %   According to Simpson~\cite{Simpson:9810} (who cites 
% %   Lewis~\cite{Lewis} as evidence), to get a sufficient amount of weakness 
% %   it is enough to 
% %   weaken one of the structures composition and identity.  This is the 
% %   idea behind Simpson's conjectures below.  However, in order to 
% %   compare fair categories with Tamsamani categories, it is necessary to
% %   go the whole length and consider fair Tamsamani categories: these are
% %   obviously functors $\grosop\to\SS$ satisfying only weak Segal 
% %   condition (i.e, preservation of fibre products over discrete 
% %   objects only up to equivalence).  It is expected that one can always
% %   stricitify one of the structures.  For example 
%   
% %   Finally there is an interesting possibility of allowing 
% %   a pseudo-discrete space (acyclic space) of objcets instead of a 
% %   truly discrete one.  This has the obvious advantage of being a 
% %   property stabel under equivalence\ldots  However, this may cause 
% %   trouble for the Segal condition: the reason the Segal condition
% %   can refer to plain fibre products instead of the `correct' notion
% %   of homotopy limits, is that the two coincide over discrete objects.
% %   With only pseudo-discrete object spaces this would no longer work, 
% %   it seems\ldots
\end{blanko}

\begin{blanko}{Fair monoids.}\label{monoid}
  A {\em fair monoid} in $\SS$ is just a fair $\SS$-category such that $O$ 
 is a singleton set $*$.  Let $\kat{S-FairMon}$ denote the full
  subcategory of $\kat{S-FairCat}$ consisting of fair monoids.  There is a
  forgetful functor $\kat{S-FairMon} \to \SS$ which sends $X = (*,A,U)$
  to its underlying space $A$.  Given two fair monoids $X=(*,A,U)$ and
  $X'=(*,A',U')$, it is easy to see that a morphism $F:X\to X'$
  is an equimorphism in $\kat{S-FairMon}$ 
  if and only if $A\to A'$ and $U\to U'$ are equimorphisms in $\SS$ 
  ---
  essential surjectivity is an empty condition.
  In view of the Segal condition and the assumption that 
  equimorphisms in $\SS$ are stable under products, this in
  turn is equivalent to requiring all components of $F$ to be
  equimorphisms, so the notion of equivalence of monoids is level-wise.
  
  To give a fair monoid in $\SS$ amounts to giving a pair of semi-monoids 
$(A,U)$ in
$(\SS,\times)$ and a semi-monoid homomorphism $U \to A$, such that $U$ is
contractible and such that the arrows
$$
U\times A\to A\leftarrow A\times U
$$ 
are equiarrows in in $\SS$.  We will refer to a fair monoid by the 
notation $(U \to C)$.

  We will use the term {\em fair monoidal category} for a fair monoid in 
  $\Cat$.
\end{blanko}

\begin{blanko}{Fair monoids in non-cartesian enriched contexts.}
%   The $\gros$-viewpoint is needed to describe the weak Segal condition
%   and to compare to Tamsamani theory.  However 
  The above
  characterisation of fair monoids, in terms of a pair of semi-monoids
  with one of them contractible and so on, might also be useful
  for the reason that it makes sense in monoidal coloured categories
  (or monoidal model categories) in which the monoidal structure is not
  the cartesian product.  For example, define a fair monoid in the
  monoidal category of chain complexes $(\kat{Ch}, \tensor, \Z)$ to
  be a pair of semi-monoids $A$ and $U$, with $U$ contractible (in the
  usual sense of homotopies of chain complexes), together with a 
  semi-monoid homomorphism $U \to A$, such that the arrows
  $$
  U\tensor A\to A\leftarrow A\tensor U
  $$ 
  are weak equivalences in $\kat{Ch}$.  I have not investigated this
  definition further.
\end{blanko}

%%%%%%%%%%%%%%%%%%%%%%%%%%%%%%%%%%%%%%%%%%%%%%%%%%
\section{Fair $2$-categories}
%%%%%%%%%%%%%%%%%%%%%%%%%%%%%%%%%%%%%%%%%%%%%%%%%%

In this section we work out the case of dimension $2$, i.e. fair
categories in $\Cat$.  In view of the canonical equivalence $\Cat
\simeq \kat{1FairCat}$ we will abuse of notation and let
$\kat{2FairCat}$ denote the category of fair categories in $\Cat$, not
in $\kat{1FairCat}$.  In the same vein, but still more abusively, by
{\em semi-$2$-category} we shall mean a semi-category enriched over
$\Cat$, i.e.~a semi-$2$-category which happens to have identity 
$2$-cells.

The main result is that a fair $2$-category is essentially the same
thing as a bicategory with strict composition law, just as one would
expect, both notions being semi-$2$-categories (in the abusive sense)
with some extra unit 
structure.   The fair category viewpoint on such a bicategory $\CC$ encodes
{\em all} possible unit structures on $\CC$ instead of favouring one of them.

\begin{blanko}{Bicategories with strict composition law.}
  A bicategory with strict composition law is just like a strict
  $2$-category, except that each object is not required to have a
  strict identity $1$-cell, but only a specified weak one.  In order
  to encourage the interpretation of bicategories as `many-object
  monoidal categories', we denote the composition law by $\tensor$ and
  compose from the left to the right.  Objects will be omitted from
  the notation whenever possible, and arrows are denoted by uppercase
  letters.  A weak identity arrow for an object $o$ is a triple
  $(I_o,\lambda,\rho)$ consisting of an arrow $I_o: o\to o$, a left
  constraint $\lambda$ and a right constraint $\rho$: these are
  invertible $2$-cells $\lambda_Y : I_o \tensor Y \isopil Y$ and
  $\rho_X : X \tensor I_o \isopil X$, natural in arrows $Y:o \to
  \cdot$ and $X:\cdot\to o$ respectively.  The left and right
  constraints are subject to the condition \vspace{2ex}
  \begin{equation}\label{Kelly}
  \begin{diagram}[w=6ex,h=4.5ex,scriptlabels,tight]
  \cdot & \rTo^X & o &\topglob{I_o \tensor Y}\lift{-2}{\Downarrow\lambda_Y}\dropglob{Y}&\cdot
  \end{diagram}
  \hspace{20pt}=\hspace{20pt}
  \begin{diagram}[w=6ex,h=4.5ex,scriptlabels,tight]
  \cdot  &\topglob{X \tensor I_o}\lift{-2}{\rho_X\Downarrow}\dropglob{X}&o& 
  \rTo^Y & \cdot
  \end{diagram}
  \end{equation}
  \vspace{2ex}
  
  A homomorphism of bicategories with strict composition law is a
  bifunctor $F:\CC\to\CC'$ which preserves the composition strictly.  The only
  non-strict part of $F$ is the comparison of identity arrows: for
  each object $o$ there is specified an invertible $2$-cell $\phi_o: I_{F(o)} 
  \isopil F(I_o)$, compatible with the respective left and right
  constraints like this:
$$
\begin{diagram}[w=6.5ex,h=5ex,tight,hug]
    I_{F(o)} \tensor F(X) && \rTo^{\lambda'_{FX}}  && F(X)  \\
   & \rdTo(1,2)_{\phi_o \tensor FX}  &    & \ruTo(1,2)_{F(\lambda_X)} & \\
   & F(I_o)\!\tensor\! F(X)  & \rLig  & F(I_o\!\tensor\! X)&
    \end{diagram}
    \hspace{24pt}
    \begin{diagram}[w=6.5ex,h=5ex,tight,hug]
    F(X)  && \lTo^{\rho'_{FX}}  &&   F(X)\tensor I_{F(o)} \\
   & \luTo(1,2)_{F(\rho_X)}  &    & \ldTo(1,2)_{FX \tensor \phi_o} & \\
   &  F(X \!\tensor\! I_o) & \rLig  & F(X) \!\tensor\! F(I_o)&
    \end{diagram}
$$
\end{blanko}

  Let $\B$ denote the category of bicategories with strict 
  composition law and bifunctors preserving the composition strictly.

\begin{prop}\label{B}
  There is an equivalence of categories
  $$
  \kat{2FairCat} \simeq \B .
  $$
\end{prop}
The equivalence is described below.  The functor $\B \to
\kat{2FairCat}$ is canonical;  the pseudo-inverse depends on a choice.
The construction relies on a couple of basic observations about identity
arrows in bicategories which are not widely known.  Further details
can be found in \cite{Kock:trivial-units}, for the one-object case.

\begin{blanko}{Identity-arrow structures.}
  Just as identity arrows in a category are uniquely determined by the
  unit axioms, in a bicategory the identity arrows are unique up to
  unique isomorphism.  To be precise, define an {\em identity arrow}
  in a semi-$2$-category to be a triple $(I_o,\lambda,\rho)$ where
  $I_o : o \to o$ is an endo-arrow of an object $o$, and $\lambda$ and
  $\rho$ are left and right constraints satisfying the axioms above.
  A {\em morphism of identity arrows} is given by a $2$-cell $I \to J$
  compatible with the left and right constraints.  Clearly the
  category of identity arrows and their morphisms is the disjoint
  union
  $$
  \Id_\CC = \coprod_{o \in \Ob(\CC)} \Id_\CC(o)
  $$
  where $\Id_\CC(o)$ is the category of identity arrows of the object $o$.
  
  The category of identity arrows has a
  natural composition law lifting the composition law on $\CC$ (and in
  particular is object-wise):
  the composition of $(I,\lambda,\rho)$ with $(I',\lambda',\rho')$
  is the composite $I\tensor I'$ equipped with left and right constraints
  \begin{diagram}[w=6ex,h=4ex,tight]
    I \tensor I' \tensor X && \rTo^{I \tensor \lambda'} & I \tensor X & 
    \rTo^{\lambda} & X   \\
     X \tensor I \tensor I' && \rTo_{\rho\tensor I'} & X \tensor I' &
     \rTo_{\rho'} & X
  \end{diagram}
  It is straightforward to check that these constraints satisfy the 
  identity arrow axiom~(\ref{Kelly}), and it also follows readily that this 
  composition law 
  is strictly associative if the original composition law is so. 
  
  Altogether, we have a semi-category $\Id_\CC$ enriched in $\Cat$, (which is 
  the
  disjoint union of the semi-monoidal categories $\Id_\CC(o)$, $o\in 
  \Ob(\CC)$).
\end{blanko}

\begin{lemma}
  The category $\Id_\CC$ of identity arrows in a bicategory $\CC$
  is equivalent to the discrete category $\Ob(\CC)$, i.e.~$\Id_\CC(o)$ is 
  contractible for each object $o$.
%   and it has a natural tensor lifting the tensor on
%   $\CC$.  If the tensor on $\CC$ is strict then so is the one on
%   $\UU$.
\end{lemma}

\begin{dem}
  Given units $(I,\lambda, \rho)$ and $(I',\lambda',\rho')$ of an 
  object $o$,
  one checks that the isomorphism
  \begin{diagram}[w=5ex,h=4ex,tight]
    I' & \rTo^{\rho^{-1}_{I'}} & I'I & \rTo^{\lambda'_I} & I
  \end{diagram}
  is compatible with the left and right constraint of $I$ and $I'$,
  hence constitutes a morphism of identity arrows, hence the category
  of identity arrows of $o$ is iso-connected.  Compatibility with the
  left and right constraints also implies that there can be at most
  one connecting arrow.
%   
%   To establish uniqueness
%   of the connecting arrows, note first that the functor given by
%   post-composing with $I$ is an equimorphism of the category $\End(o)$.
%   Indeed, the right constraint provides a comparison isomorphisms to
%   the identity functor.  Now suppose there are two connecting
%   morphisms $I'\topile I$.  After post-composition with $I$ these two
%   $2$-cells will coincide, as a consequence of compatibility with left
%   constraints, and because the involved arrows are invertible.  But
%   since post-composition with $I$ is an equivalence, it is a bijection
%   on the level of $2$-cells, hence already the two original morphisms
%   coincide.
\end{dem}

\begin{blanko}{From bicategories to fair $2$-categories.}
  With the above observations it is easy to define the functor $\B \to
  \kat{2FairCat}$.  In analogy with the fair nerve described 
  in~\ref{fairnerve}, given a bicategory $\CC$ with strict
  composition we define a $\grosop$-diagram $X=(O,A,U)$
  in $\Cat$ like this:
  \begin{eqnarray*}
    \grosop & \longrightarrow & \Cat  \\
    \intextB & \longmapsto & O \df \Ob \CC \\
    \intextBB & \longmapsto & A \df \coprod_{x,y\in\Ob \CC} \Hom(x,y) \\
    \intextBbB & \longmapsto & U \df \coprod_{x\in\Ob\CC} \Id_\CC(x)
  \end{eqnarray*}
%   Here the category $\operatorname{Identities}(x)$ is defined as above.
%  It is  contractible for each $x$ and has a natural tensor, as required.
  The map $U \to A$ is the forgetful functor sending an identity
  triple $(I,\lambda,\rho)$ to its supporting arrow.  The rest of the
  $\grosop$-diagram is defined by invoking the strict Segal condition
  and observing that $U \to A$ is a $\tensor$-functor.
 
  By definition of the morphisms in $\Id_\CC$, the collection of all left
  constraints at $o$ assembles into an
  invertible natural transformation
  \vspace{2ex}
  \begin{diagram}[w=6ex,h=4.5ex,tight]
  \Id_\CC(o) \times \Hom_\CC(o,\cdot) &&& 
  \topglob{\tensor}\lift{-2.2}{\Updownarrow \lambda} 
  \dropglob{\text{\footnotesize proj}}  && \Hom_\CC(o,\cdot) \ .
  \end{diagram}
  \vspace{1ex}

  \noindent The component on a unit $(I,\lambda,\rho)$ and an arrow
  $X$ is nothing but $\lambda_X: I \tensor X \to X$.  Since
  $\Id_\CC(o)$ is contractible, the projection map is an equimorphism,
  and hence the composition map is too.  Summing over all the objects 
  we get an
  invertible $2$-cell 
  
  \vspace{1.4ex}
  \begin{diagram}[w=5ex,h=4.5ex,tight]
  U \times_O A && \topglob{\tensor}\lift{-2.4}{\Updownarrow \lambda} 
  \dropglob{\text{\footnotesize proj}}  & A
  \end{diagram}
  \vspace{1ex}
  
  \noindent
  hence the composition map is an equimorphism as required.
    Similarly for the right constraints.

    The same argument works for any vertical arrow in
    $\gros$: it is the dot-sum ($\dotsum$) of identity arrows and the case just 
    treated.
     
    \bigskip
    
    One can check that this construction is functorial
    (cf.~\cite{Kock:trivial-units}): the unit part $\phi_o$ of a
    bifunctor $F: \CC\to\DD$ amounts precisely to a lift of $F$ to
    $\Id_\CC\to\Id_\DD$, which in turn is equivalent to extending the
    natural transformation between the corresponding
    $\Dmono\op$-diagrams to a natural transformation of
    $\grosop$-diagrams.  This finishes the construction of the functor
    $\B\to \kat{2FairCat}$.
\end{blanko}

\begin{blanko}{From fair $2$-categories to bicategories.}
  Given a fair category $(O,A,U)$ in $\Cat$, the $A$-part already 
  constitutes a semi-$2$-category $A$, whose composition we denote
  by $\tensor$.  
  It remains to use
  the $U$-part of the diagram to provide weak identity arrows for the
  semi-$2$-category $A$.  Note that $U$ also constitutes a 
  semi-category in $\Cat$; again we denote the composition law by 
  $\tensor$; then the functor $U \to A$ is a $\tensor$-functor.  The
  equivalence $U \isopil O$ has a pseudo-section (since we are just
  talking plain categories), i.e., for each element $o\in O$ we can
  pick an object $I_o$ in $U$, a chosen weak unit. 
  The functor we are constructing depends on this choice, but different 
  choices will yield canonically isomorphic results.
  
  Since $U$ is 
  contractible,
  
  \medskip
  
  (S1) \ \label{S1S2}
  we have an isomorphism $\alpha: I_o \tensor I_o \isopil I_o$,
  
  \medskip
  
  \noindent and since the composition functors  
  $ U \times_O A \stackrel{\tensor}{\rTo} A \stackrel{\tensor}{\lTo} A \times_O U$
  are equimorphisms of categories, we see that 
  
  \medskip
  
  (S2) 
  composition with $I_o$ defines equimorphisms
  $A(o,\cdot) \!\to\! A(o,\cdot)$ and $A(\cdot,o) \!\to\! A(\cdot,o)$.
  
  \medskip
  
  \noindent
  This means that for each arrow $X:o\to \cdot$ there is a bijection
  $$
  \operatorname{2Cell}_\CC (I\!\tensor\! X,X) \simeq 
  \operatorname{2Cell}_\CC( I\!\tensor \!I\!\tensor\! X, I\!\tensor\! X) .
  $$
  On the right-hand side we have the canonical $2$-cell $\alpha
  \tensor X$, so take the inverse image and call it $\lambda_X$, the
  required left constraint.  The right constraint is constructed
  analogously.  It is easy to see that these constraints are natural,
  and one can also check that they satisfy the coherence
  condition~(\ref{Kelly}).  Both these claims follow from the fact that
  $\lambda$ and $\rho$ are defined in terms of the isomorphism
  $\alpha$ in $U$, which is automatically coherent since it lives in a
  contractible category.  (See \cite{Kock:trivial-units} for further 
  details.)
\end{blanko}

\begin{BM}\label{Saa}
  Conditions (S1) and (S2) provide in fact a useful
  definition of unit.  This viewpoint goes back to 
  Saavedra~\cite{Saavedra} in the case of monoidal categories; it is
  exploited further in \cite{Kock:trivial-units}.  The relevance of
  this viewpoint in higher dimensions was first suggested by
  Simpson~\cite{Simpson:9810}.  The basic $3$-dimensional theory
  is worked out in Joyal-Kock~\cite{Joyal-Kock:coherence}.
\end{BM}

\begin{blanko}{Non-strict bifunctors.}
  \label{weak-morphisms} Even for strict $2$-categories it is
  sometimes necessary to consider weaker notions of $2$-functors
  (strong or lax), which do not respect composition strictly, but only
  up to specified comparison arrows (isomorphisms or arbitrary
  morphisms), cf.~B\'enabou~\cite{Benabou:bicat}.  Such $2$-functors
  can also be captured in the setting of fair $\Cat$-categories: they
  correspond precisely to weaker natural transformations (strong or
  lax) between $\grosop$-diagrams $u:F\Rightarrow G$, satisfying a
  strict `Segal condition': $u_{m\smalldotsum n} = u_m \times_{u_O} u_n$.
  (Note that any weak natural transformation satisfies this condition
  up to isomorphism, and is equivalent to one satisfying it strictly.)
\end{blanko}

%%%%%%%%%%%%%%%%%%%%%%%%%%%%%%%%%%%%%%%%%%%%%%%%%%
\section{Fair $3$-categories}
%%%%%%%%%%%%%%%%%%%%%%%%%%%%%%%%%%%%%%%%%%%%%%%%%%

\label{dim3}

For more details on the results and constructions in this section,
see Joyal-Kock~\cite{Joyal-Kock:coherence}.

\bigskip

Gordon, Power, and Street~\cite{Gordon-Power-Street} introduced the
notion of tri\-categories, designed to be the weakest possible
definition of $3$-dimensional category.  In analogy with the
$2$-dimensional case, fair $3$-categories should correspond to
tri\-categories with strict composition laws.  Carrying out the
comparison in the style of the previous section would seem to be somewhat
involved, though.
% The main difficulty is that the values of the
% $\grosop$-diagrams no longer possess strict identity arrows, so for
% example given $U \times_O A \to A$, it is not possible to fix the
% first argument and obtain $A \to A$.  That construction implicitly
% involves identity arrows in $U$.
What we undertake here is a comparison in the semi-strict situation of
fair categories in $\kat{2Cat}$ (interpreted as a semi-strict sort of
fair $3$-categories in the canonical way, cf.~\ref{fair-n}), which we
compare to locally strict tri\-categories with strict composition laws.
Furthermore for the sake of clarity, we restrict to the one-object
version.

\begin{blanko}{Strictly monoidal strict $2$-categories with weak 
  units.}\label{GPS-strict}
  Starting with the definition in \cite{Gordon-Power-Street}, we 
  specialise to the
  locally strict case (i.e.~every hom bicategory is a strict 
  $2$-category). Then we restrict to the one-object case, and finally 
  we require the composition law to be strict.  The resulting notion 
  is a strict $2$-category $\CC$ with a strict tensor product with a weak 
  unit.  The notion of weak unit in this situation is a specified triple
  $(I,\lambda,\rho)$ where $I$ is an object in $\CC$,
$\lambda$ is a natural family of equimorphisms (i.e.~admitting a
pseudo-inverse) $\lambda_X :I \tensor X \to X$, and $\rho$ is a 
natural family of equimorphisms $\rho_X: X \tensor I \to
X$.  The naturality is only up to isomorphism, so  
specifying the lambdas and rhos involves specifying certain invertible
$2$-cells.  Finally there is specified a natural family of invertible 
$2$-cells (a modification)
$$
\cel K :\rho_X \tensor Y \Rightarrow X \tensor \lambda_Y
$$
subject to two coherence axioms called left and right normalised 
$4$-cocycle conditions (\cite{Gordon-Power-Street}, 
p.11 and p.12), which we do not reproduce here.
\end{blanko}

\begin{blanko}{Fair monoids in $\kat{2Cat}$.}
  Let $(U\to C)$ denote a {\em fair monoid} in $\kat{2Cat}$, 
  cf.~\ref{monoid}.  This means that
  
  --- $U$ and $ C$ are strict semi-monoidal $2$-categories, 
  
  --- there is a strict semi-monoidal $2$-functor $U \to C$,
  
  ---  $U$ is contractible, and the multiplication maps
  $$
  U \times C \to C \leftarrow C \times U
  $$
are equi-$2$-functors.  
\end{blanko}

\begin{prop}\label{corr}
    Fair monoids in $\kat{2Cat}$ correspond to strict 
    $2$-categories with strict tensor product and weak unit.
\end{prop}
The proof is comprised by the following two subsections.

\begin{blanko}{From fair monoid in $\kat{2Cat}$ to monoidal $2$-category.}
  Given a fair monoid $(U \to C)$, since $U$ is contractible, for each
  object $I$, there exists an arrow $\alpha : I \tensor I \isopil I$
  in $U$; this arrow is an equiarrow, and any two such are uniquely
  isomorphic.  We fix $I$ and $\alpha$, and we use the same symbols
  for their images in $C$.  It also follows from contractibility of 
  $U$ that there is an associator modification $\cel A: I \tensor 
  \alpha \isopil \alpha \tensor I$ which satisfies the pentagon 
  equation.

  As in the $2$-dimensional case, the fact that the multiplication
  map $U \times C \to C$ is an equi-$2$-functor
  means that the same is true for the $2$-functor `tensoring with
  $I$ from the left' (and similarly with tensoring from the right).
  Therefore, for each object 
  $X$ in $C$ there is an
  equimorphism of categories
$$
\Hom( I\tensor X , X) \isopil 
\Hom(I \tensor I \tensor X , I \tensor X)   .
$$
In the second category we have the canonical object $\alpha \tensor
X$.  Hence there exists a pseudo pre-image $\lambda_X$, together with
an isomorphism $\cel L : I \tensor \lambda_X \Rightarrow \alpha \tensor
X$.  For chosen $\lambda_X$ and $\cel L$, there is a unique way to
assemble the lambdas into a natural transformation (this involves
specifying some $2$-cells), in such a way that $\cel L$ becomes 
natural in $X$.  The pair $(\lambda,\cel L)$ is not unique, but
any two such are uniquely isomorphic.

Similarly there is a natural transformation $\rho$ with components
$\rho_X : X \tensor I  \to X$ equipped with a natural modification
$\cel R : \rho_X \tensor I \Rightarrow X \tensor \alpha$ (and this data is
unique up to unique isomorphism).

The lambdas and rhos are the required left and right constraints.
Finally, using $\cel L$, $\cel R$, and $\cel A$ one can construct a
modification $\cel K : \rho_X \tensor Y \Rightarrow X \tensor
\lambda_Y$ which satisfies the normalised $4$-cocycle conditions of
\cite{Gordon-Power-Street} --- this is a consequence of the pentagon
equation for $\cel A$, and hence ultimately a consequence of
contractibility of $U$.  See \cite{Joyal-Kock:coherence} for details.

\bigskip

There were many choices involved in the construction: first the choice
of $I$ and $\alpha$, then the choices of $\lambda$ and $\rho$, 
together with the auxiliary $2$-cells $\cel L$ and $\cel R$.  But all
these choices lead to equivalent monoidal $2$-categories.

% This is the definition of monoidal $2$-category that can be extracted
% from the notion of fair monoid in $\kat{2Cat}$.  The modifications $L$
% and $R$ express the failure of the actions $\lambda$ and $\rho$ to be
% strictly associative.  One then verifies that these two `associators'
% satisfy the relevant pentagon conditions (which involve also the
% `original' associator $A$), and also that there is a canonical
% modification expressing the up-to-equivalence bimodule condition for
% compatibility between the two actions.  
% All this expresses what a monoidal $2$-category is according to the
% theory of  fair categories.  Here, a very strict version of monoidal
% $2$-category is explained, but you can weaken the underlying
% $2$-categories to be bicategories, and you can allow the tensors to be
% weak, with associators satisfying pentagon equation up to Stasheff's $K3$
% condition.  

\end{blanko}

\begin{blanko}{From monoidal $2$-category to fair monoid in $\kat{2Cat}$.}
  Starting from a monoidal $2$-category $\CC$ like specified in
  \ref{GPS-strict}, ideally one would construct the $2$-category of
  all unit structures on $\CC$.  The objects are quadruples
  $(I,\lambda,\rho,\cel K)$ like in \ref{GPS-strict}, and arrows and
  $2$-cells are defined to come equipped with compatibility data with
  respect to these structures.  This category is contractible as
  required, but unfortunately its tensor product is not strict as
  required in order to play the r\^ole of $U$ in a fair monoid
  $(U \to C)$.
  
  An alternative construction is used.  Define $\alpha: I \tensor I
  \isopil I$ by taking $\alpha\df\lambda_I$.  (We
  could equally well have used $\rho_I$: there is a canonical 
  modification $\lambda_I \Leftrightarrow \rho_I$ constructed from 
  $\cel K$ and the naturality data specified with $\lambda$ and $\rho$.) 
  Now
  the $2$-functors $\CC\to\CC$ defined by multiplying with $I$ from
  the left or the right are equi-$2$-functors.  Hence the pair
  $(I,\alpha)$ satisfies the $2$-dimensional version of (S1) and (S2)
  on page \pageref{S1S2}.
  It is shown in \cite{Joyal-Kock:coherence} that $\alpha:I\tensor I 
  \to I$ is associative up to a canonical associator modification 
  $\cel A : I \tensor \alpha \isopil \alpha \tensor I$,
  which satisfies the pentagon equation.  Hence $(I,\alpha)$ is a 
  weak semi-monoid. 
  The set of such
  pairs $(I,\alpha)$ is the object set of a contractible $2$-category,
  the arrows being weak-semi-monoid equimorphisms, and the $2$-cells 
  being invertible $2$-cells in $\CC$.  Again, however, the tensor
  structure on this $2$-category is not strict.  But in this case
  there is a strict alternative: it is the free strictly semi-monoidal
  $2$-category whose objects are the positive powers of $I$, whose
  generating arrows are $\alpha$ and a fixed right adjoint $\beta$, and whose 
  generating $2$-cells are $\cel A$ together with unit and counit for 
  $\alpha\isleftadjointto \beta$.  This
  $2$-category is contractible and does the job as $U$.  The canonical
  $2$-functor consisting in interpreting this category in $\CC$
  respects the tensor product.  Now we have got semi-monoidal
  $2$-categories $U$ and $C$ and a semi-monoidal $2$-functor between
  them, so the rest of the $\grosop$-diagram is determined by the
  strict Segal condition.
  
  Just like in the $2$-dimensional case, the left and right
  constraints provide natural $2$-cells between
  multiplication-with-a-unit and projection, showing that $U\times A
  \to A \leftarrow A \times U$ are equivalences.  Hence the
  $\grosop$-diagram is a fair monoid.
\end{blanko}

% \pagebreak
%%%%%%%%%%%%%%%%%%%%%%%%%%%%%%%%%%%%%%%%%%%%%%%%%%
\section{$n$-groupoids, homotopy $n$-types, and Simpson's conjecture}
%%%%%%%%%%%%%%%%%%%%%%%%%%%%%%%%%%%%%%%%%%%%%%%%%%

\label{Sec:groupoids}

\begin{blanko}{Motivation.}
  Important motivation for higher category theory is the desire of
  giving a completely algebraic account of homotopy theory
  (cf.~Grothendieck~\cite{Grothendieck:stacks}).  In particular, every
  topological space should have associated a {\em higher fundamental
  groupoid}.  This should be an $\infty$-groupoid, but truncated
  homotopy $n$-types should be described by $n$-groupoids.  Ideally
  this description should be a pair of (weakly) adjoint functors
$$
R : \infty\Grpd \pile{ \rTo \\ \lTo} \Top : \Pi
$$
where $\Pi$ is the fundamental groupoid mentioned above, and $R$ is a
geometric realisation functor, and for each $n\geq 0$
this pair of functors should induce an
equivalence between the homotopy category of $n$-groupoids and that of
$n$-truncated topological spaces.
\end{blanko}

\begin{blanko}{Completely strict $n$-groupoids.}
  The strictest possible definition of $n$-groupoid is as a strict
  $n$-category such that every $k$-cell is strictly invertible, for $k\leq
  n$.  It was shown quite early
  (Brown-Higgins~\cite{Brown-Higgins:infty-groupoids}) that such strict
  $n$-groupoids cannot realise all homotopy $n$-types.  The problem occurs
  already with the $3$-type of $S^2$:  this space has a non-trivial
  Whitehead operation $\pi_2\tensor \pi_2 \to \pi_3$, but a version of the
  Eckmann-Hilton argument shows that every strict $n$-groupoid gives trivial
  Whitehead brackets.  (See Simpson~\cite{Simpson:9810} for a detailed
  account of these arguments.)  
  
  A slightly weaker notion of strict $n$-groupoid is obtained by requiring
  the $k$-cells to be invertible only up to a $(k+1)$-cell, which in turn
  should be weakly invertible in the same sense, and so on, up to the
  $n$-cells which should be strictly invertible.  The following formulation
  of this idea can be interpreted in either strict $n$-categories,
  Tamsamani $n$-categories, or fair $n$-categories:
\end{blanko}
  
\begin{blanko}{$n$-groupoids.}\label{n-groupoids}
  An $n$-category $X$ is called an {\em
  $n$-groupoid} if the category $\tau_0 X$ is a groupoid, and if for each
  pair of objects $x,y$, the $(n-1)$-category $X(x,y)$ is an
  $(n-1)$-groupoid.
\end{blanko}

\begin{blanko}{Kapranov-Voevodsky $n$-groupoids.}
  If in the above definition, `$n$-category' is taken to mean `strict
  $n$-category' then the $n$-groupoid notion is that of
  Kapranov-Voevod\-sky~\cite{Kapranov-Voevodsky}.  Since the composition laws
  as well as the identity cells are strict, the Eckmann-Hilton argument
  still applies, showing that the standard geometric realisation of such an
  $n$-groupoid will have trivial Whitehead brackets.  Kapranov and
  Voevodsky constructed a new realisation functor with certain properties
  (preservation of homotopy groups), and
  claimed that with this, their $n$-groupoids could realise all homotopy
  $n$-types.  Simpson~\cite{Simpson:9810} showed instead that {\em any}
  realisation functor satisfying those properties will yield only homotopy
  types with trivial Whitehead brackets, and concluded that there must be
  an error in either \cite{Kapranov-Voevodsky} or in \cite{Simpson:9810}.
\end{blanko}

\begin{blanko}{Tamsamani $n$-groupoids.}\label{Tams-thm}
  If in Definition~\ref{n-groupoids}, `$n$-category' is taken in the
  sense of Tamsamani, the situation is different:
  Tamsamani~\cite{Tamsamani:thesis} has constructed the fundamental
  $n$-groupoid of an ($n$-truncated) topological space, and a left
  adjoint realisation functor, and he has shown that this adjunction
  induces an equivalence of homotopy categories between $n$-groupoids
  in his sense and $n$-truncated topological spaces.  (The existence
  of this construction is one of the main advantages of his
  definition over many of the other existing definitions of weak
  higher categories.)
\end{blanko}

\begin{blanko}{Simpsons conjecture(s).}\label{Simpson-conj}
  Although Tamsamani's theorem shows that weakening the notion of
  composition suffices to capture all homotopy $n$-types, it is a
  puzzling idea that the real problem with strictness are the strict
  identity arrows: the Whitehead operations are trivial because of the
  Eckmann-Hilton argument, which in turn relies crucially on
  strict identities.  A detailed analysis of these issues led
  Simpson~\cite{Simpson:9810} to formulate the following
  conjecture(s):
  
  {\em There exists a notion of strict $n$-groupoid with weak identity 
  arrows and a notion of geometric realisation such that every 
  homotopy $n$-type appears in this way.  
  
  In the other direction there should be a fundamental $n$-groupoid
  (with strict composition and weak identities) associated to every
  topological space, and these two functors should induce an
  equivalence between the homotopy categories of $n$-groupoids with
  weak identities on one side and $n$-truncated topological spaces on
  the other side.
  
  More generally, the homotopy theory of strict $n$-categories with weak 
  identities should be equivalent to the homotopy theory of Tamsamani 
  $n$-categories.}
  
  An ad hoc definition of weak identities was sketched in the preprint
  (based on (S1) and (S2) from page \pageref{S1S2}), 
  but it was acknowledged that it might not be
  the correct definition to turn the conjecture true, and in fact, the
  details of this definition were never worked out.  Simpson's
  conjecture was one starting point for the present work, and the
  notion of fair category emerged gradually from an attempt to
  understand his ideas.
  
  The conjecture in its strong form has startling consequences,
  defying all trends in higher category theory: every 
  weak higher category should be equivalent to one with strict
  composition!
  
  As an annex to Simpson's conjecture, and at the same time a
  concretisation of the objects of its assertion, I want to
  propose that the notion of fair $n$-groupoid is appropriate for
  fulfilling the conjecture.  For emphasis, here is the definition:
\end{blanko}

\begin{blanko}{Fair groupoids.}\label{fairgroupoids}
  A fair $n$-category $X:\grosop\to \kat{(n-1)FairCat}$ is called a
  {\em fair $n$-groupoid} if the category $\tau_0 X$ is a groupoid and
  if for each pair of objects $x,y$, the fair $(n-1)$-category
  $X(x,y)$ is a fair $(n-1)$-groupoid.
  
  (Note that this definition does not explicitly refer to the weak 
  identity arrows, but that it relies on \ref{1FCat} as base for the 
  induction.)
\end{blanko}

%%%%%%%%%%%%%%%%%%%%%%%%%%%%%%%%%%%%%%%%%%%%%%%%%%
\subsection*{Simpson's conjecture in dimension $3$}
%%%%%%%%%%%%%%%%%%%%%%%%%%%%%%%%%%%%%%%%%%%%%%%%%%

The crucial test for the conjecture is dimension $3$.  If $*$ is an
object in a (strict or weak) $3$-groupoid $G$, and $I$ is a (weak)
identity arrow of $*$, then $\End(I)$ is a braided monoidal category
(in fact a braided categorical group).  The braiding corresponds to 
the
Whitehead operation $\pi_2 \tensor \pi_2 \to \pi_3$ (relative to the 
base point corresponding to $*$).
This Whitehead bracket is zero if and only if the
braiding collapses to a symmetry.  So the failure of strict
$3$-groupoids to realise all $3$-types is equivalent to the fact that
if $G$ is strict then the braiding of $\End(I)$ is a symmetry (in
fact, is equivalent to a commutativity).
There is no essential generality lost in treating only the
one-object case, so we consider $1$-object $3$-groupoids, 
which is the same thing as monoidal $2$-groupoids such that every 
object has a weak tensor inverse;  let $I$ denote the (weak) unit
for the tensor product.  Simpson's conjecture in this case says that 
such monoidal $2$-groupoids with strict composition laws and strict
tensor 
can realise all pointed $3$-types. 
Let us further restrict 
attention to
simply connected $3$-types: this corresponds to having a contractible
space of objects, so all objects are equivalent to $I$.
The following result is a form of Simpson's conjecture in this case:

\begin{satz}\label{JK}
  (Joyal-Kock~\cite{Joyal-Kock:traintracks}.)  Every braided monoidal
  category arises as $\End(I)$, where $I$ is a weak unit in an
  otherwise completely strict monoidal $2$-category (as treated in 
  Proposition~\ref{corr}).
\end{satz}

$1$-connected homotopy $3$-types correspond to braided categorical
groups.  Under the correspondence of the theorem, these correspond to
strict $2$-groupoids with invertible tensor product and weak units
(which in turn are one-object fair $3$-groupoids in the sense of 
\ref{fairgroupoids}).
Hence we get the following version of Simpson's conjecture in
dimension $3$:

\begin{cor}\label{JKcor}
  (Cf.~\cite{Joyal-Kock:traintracks}.)
  Strict $2$-groupoids with invertible tensor product and weak units
  can model all $1$-connected homotopy $3$-types.
\end{cor}

The idea of the proof of the theorem is this: the braiding from $g\circ f$
to $f \circ g$ relies on the arrow $\alpha : I \tensor I \to I$, 
cf.~the proof of \ref{corr}, and
a chosen quasi-inverse.  In graphical notation (reading the string 
diagrams from the bottom to the top), we picture $\alpha$ as
\raisebox{-4pt}{
  \begin{texdraw}
  \drawdim pt \textref h:C v:C \setunitscale 3
  \arrowheadsize l:0.32 w:0.24 \arrowheadtype t:F
  \move (-2.5 0) \move (1.5 0)
\move (0 0) \bsegment
  
  \lvec (0 -1.2)
  \clvec (0 -3)(1 -3)(1 -4.8) \lvec (1 -6)
  \move (0 -1.2)
  \clvec (0 -3)(-1 -3)(-1 -4.8) \lvec (-1 -6)
  \esegment
\end{texdraw}
}
Then the braiding is this:

\newcommand{\vigespor}{%
	\bsegment
	\move (0 0)
	\lvec (0 3)
	\clvec (0 5)(1 5)(1 7)
	\lvec (1 8)
	\clvec (1 10)(0 10)(0 12)
	\move (0 3)
	\clvec (0 5)(-1 5)(-1 7)
	\lvec (-1 8)
	\clvec (-1 10)(0 10)(0 12)
	\lvec (0 15)
	\savepos (0 15)(*ex *ey)
	\esegment
	\move (*ex *ey)
}

\newcommand{\rightlabeldot}[1]{%
 \bsegment
 \fcir f:0 r:0.3 \htext (0.8 0){{\footnotesize ${#1}$}}
 \esegment
}
\newcommand{\leftlabeldot}[1]{%
 \bsegment
 \fcir f:0 r:0.3 \htext (-0.8 0){{\footnotesize ${#1}$}}
 \esegment
}

\begin{center}
  \begin{texdraw}
  
  \drawdim pt \textref h:C v:C
  \setunitscale 8

  \arrowheadsize l:0.32 w:0.24
  \arrowheadtype t:F
  
%-----------------------
\move (0 0) \bsegment

%   \move (-1.5 -6) \bsegment \linewd 0.03
%     \move (0 3) \avec (0 0.2)
%     \move (0.2 0) \avec (3 0)
%   \esegment
  
  \move (0 0) \bsegment
  \lvec (0 15)
  \move (0 11) \rightlabeldot{g}
  \move (0 4) \leftlabeldot{f}
  \esegment
  
\esegment
%-----------------------
  \htext (3.75 14.1) {$\Rightarrow$}
%-----------------------
\move (7.5 0) \bsegment

%   \move (-1.5 -6) \bsegment \linewd 0.03
%     \move (0 3) \avec (0.0 0.5)
%     \move (0.1 0.3) \avec (0.9 1.1)
%     \move (1.1 0.9) \avec (0.3 0.1)
%     \move (0.5 0.0) \avec (3 0)
%   \esegment
  
  \move (0 0) \bsegment
    \vigespor
    \move (0 13) \rightlabeldot{g}
    \move (0 2) \leftlabeldot{f}
  \esegment
  
\esegment
%-----------------------
  \htext (11.25 14) {$\Rightarrow$}
%-----------------------
\move (15 0) \bsegment

%   \move (-1.5 -6) \bsegment \linewd 0.03
%     \move (0.0 3.0) \avec (0.9 2.1)
%     \move (1.0 2.0) \avec (1.0 1.0)
%     \move (0.9 0.9) \avec (0 0)
%     \move (0.4 0) \avec (3 0)
%   \esegment
  
  \move (0 0) \bsegment
  \vigespor
  \move (1 6.5) \rightlabeldot{f}
  \move (0 13) \leftlabeldot{g}
  \esegment
  
\esegment
%-----------------------
  \htext (18.75 14) {$\Rightarrow$}
%-----------------------
\move (22.5 0) \bsegment

%   \move (-1.5 -6) \bsegment \linewd 0.03
%     \move (0 3) \avec (0.9 2.1)
%     \move (1.0 2.0) \avec (1.0 1.1)
%     \move (1.1 1.0) \avec (2.0 1.0)
%     \move (2.1 0.9) \avec (3 0)
%   \esegment
  
  \move (0 0) \bsegment
  \vigespor
  \move (1 6.5) \rightlabeldot{f}
  \move (-1 8.5) \leftlabeldot{g}
  \esegment
  
\esegment
%-----------------------
  \htext (26.25 14) {$\Rightarrow$}
%-----------------------
\move (30 0) \bsegment

%   \move (-1.5 -6) \bsegment \linewd 0.03
%     \move (0 3) \avec (0.9 2.1)
%     \move (1.0 2.0) \avec (1.9 2.0)
%     \move (2.0 1.9) \avec (2.0 1.0)
%     \move (2.1 0.9) \avec (3 0)
%   \esegment
  
  \move (0 0) \bsegment
  \vigespor
  \move (1 8.5) \rightlabeldot{f}
  \move (-1 6.5) \leftlabeldot{g}
  \esegment
  
\esegment
%-----------------------
  \htext (33.75 14) {$\Rightarrow$}
%-----------------------
\move (37.5 0) \bsegment

%   \move (-1.5 -6) \bsegment \linewd 0.03
%     \move (0.0 3.0) \avec (0.9 2.1)
%     \move (1.0 2.0) \avec (2.0 2.0)
%     \move (2.1 2.1) \avec (3 3)
%     \move (3 2.6) \avec (3.0 0.0)
%   \esegment
  
  \move (0 0) \bsegment
  \vigespor
  \move (0 13) \rightlabeldot{f}
  \move (-1 6.5) \leftlabeldot{g}
  \esegment
  
\esegment
%-----------------------
  \htext (41.25 14) {$\Rightarrow$}
%-----------------------
\move (45 0) \bsegment

%   \move (-1.5 -6) \bsegment \linewd 0.03
%     \move (0 3) \avec (2.5 3)
%     \move (2.7 2.9) \avec (1.9 2.1)
%     \move (2.1 1.9) \avec (2.9 2.7)
%     \move (3 2.5) \avec (3 0)
%   \esegment
  
  \move (0 0) \bsegment
  \vigespor
  \move (0 13) \rightlabeldot{f}
  \move (0 2) \leftlabeldot{g}
  \esegment
  
\esegment
%-----------------------
  \htext (48.75 14.2) {$\Rightarrow$}
%-----------------------
\move (52.5 0) \bsegment

% \move (-1.5 -6) \bsegment \linewd 0.03
%   \move (0 3) \avec (2.8 3)
%   \move (3 2.8) \avec (3 0)
% \esegment
  
  \move (0 0) \bsegment
  \lvec (0 15)
  \move (0 11) \leftlabeldot{f}
  \move (0 4) \rightlabeldot{g}
  \esegment
  
\esegment

\end{texdraw}\end{center}

The picture suggests that going left past each other is not 
the same as going right past each other, and hence the braiding should 
not be a symmetry.  The proof that every braided monoidal category 
arises as such an $\End(I)$ consists in taking 
these diagrams seriously, relating them (in an up-to-homotopy sense)
to the geometry of labelled configuration spaces.

\bigskip

Theorem~\ref{JK} provides a direct check of the first non-trivial case
of Simpson's conjecture, and by expressing the arguments in terms of
familiar mathematical objects like braided monoidal categories it also
provides good intuitive insight to the problem.  However it comes
short in providing the notions and tools necessary for generalisation
to higher dimension: at present, the appropriate notions of geometric
realisation of fair categories have not been worked out, and there is
no general construction of fundamental fair groupoid to provide a
functor in the other direction (but see however the first example in
the next section).

% \pagebreak

%%%%%%%%%%%%%%%%%%%%%%%%%%%%%%%%%%%%%%%%%%%%%%%%%%
\section{A couple of examples}
%%%%%%%%%%%%%%%%%%%%%%%%%%%%%%%%%%%%%%%%%%%%%%%%%%
\label{Sec:examples}

The first example, somewhat detailed, concerns Moore path spaces.  It
is a fair category in $\Top$.  The second and third examples are more
succinct: the second example is about cofibrant objects in a monoidal
model category; the final example is a fair Tamsamani $2$-category of
oriented cobordisms.

\newcommand{\R}{\mathbb{R}}
%%%%%%%%%%%%%%%%%%%%%%%%%%%%%%%%%%%%%%%%%%%%%%%%%%
\subsection*{Moore path spaces}
%%%%%%%%%%%%%%%%%%%%%%%%%%%%%%%%%%%%%%%%%%%%%%%%%%

This first example is a fair category in $\SS=\Top$, where for simplicity 
we take the equimorphisms to be the homotopy equivalences.

\begin{blanko}{Moore paths.} 
  Let $X$ be a topological space.  A {\em Moore path}
  in $X$ is a continuous
  map $[0,r]\to X$, where $r>0$.  When composing paths by concatenation
  of intervals (see below), the length of the domain interval
  increases; since there is no reparametrisation involved, the
  composition is strictly associative.  In this way, taking the set
  of points of $X$ as object set $O$, and taking as arrows the space $A$
  of all Moore paths (with the compact-open topology), we have got a
  topological semi-category $O \topileback A$.
  
  Since the length of the domain interval increases strictly in every
  composition, there can be no strict unit paths.  One could of course
  just allow the zero-length interval as domain for a path, then these
  paths would be strict units.  But note that for each point $x\in X$
  the space of Moore loops at $x$ of length $r>0$ is isomorphic to the
  classical loop space $\Omega_{x}$, while for $r=0$ we get just a point.
  Furthermore, since we want to think of the path space as a groupoid, the
  null-homotopic loops should actually also be considered unit
  paths --- weak unit paths, that is.
\end{blanko}
  
\begin{blanko}{A topological fair category of Moore paths.}
  We would like to consider all null-homotopic loops at a point $x\in
  X$ as weak identity arrows of $x$.  However, the space of
  null-homotopic loops at $x$ is not contractible, as required by the
  colour axiom \ref{colour-axiom}.  What is missing is of course to 
  specify in precisely which sense such a null-homotopic loop
  serves as weak unit: as $U(x)$ we need to take the space of all 
  null-homotopic loops based at
  $x$, {\em together with a null-homotopy}.  Put $U \df \coprod_{x\in 
  O} U(x)$.
\end{blanko}

\begin{prop}
  The triple $(O,A,U)$ of points, Moore paths, and 
  null-homotopic-Moore-loops-with-given-null-homotopy, as defined above 
  constitutes a fair category 
  in $\Top$.
\end{prop}

\begin{dem}
  We first give an explicit description of the semi-categories $A$ and $U$
  and the semi-functor $U \to A$;
  then we show that $U(x)$ is contractible for each $x$, and finally
  that the maps $U \times_O A \topile  A \topileback A \times_O U$
  are equivalences.  These verifications are pretty straightforward,
  but it is instructive to see in detail how the homotopies built 
  into the individual weak units assemble into the global equivalences
  required by the axioms of fair category.
  \nomoreqed
\end{dem}
  
\begin{blanko}{Specific description of the semi-categories $A$ and $U$.}
  Given two Moore paths $\gamma_1 : [0,r_1] \to X$ and $\gamma_2 :
  [0,r_2] \to X$ with $\gamma_1(r_1)=\gamma_2(0)$, their composite
  is the Moore path $\gamma_1 \tensor \gamma_2$ defined by
  \begin{eqnarray*}
    \gamma_1 \tensor \gamma_2: [0,r_1+r_2] & \longrightarrow & X  \\
    t & \longmapsto & 
    \begin{cases} 
      \gamma_1(t) & \text{ for } t\leq r_1 \\ 
      \gamma_2(t-r_1) & \text{ for } t\geq r_1 .
      \end{cases}
  \end{eqnarray*}
  It is clear that this composition law is strictly associative.
  
  Anticipating the definition of $U$, define 
  the interval $[0,r]$ to be the segment of the 
  $(y\!=\!1)$-line in $\R^2$ starting at $(\begin{smallmatrix}0\\ 
  1\end{smallmatrix})$ and ending at $(\begin{smallmatrix}r\\ 1
  \end{smallmatrix})$:
  \begin{center}
    \begin{texdraw}
      \setunitscale 20
      \linewd 0.02
      
      % coordinate axes
      \move (-0.5 0) \lvec (6 0)
      \move (0 -0.5) \lvec (0 2.5)

      \htext (-0.5 1) {\footnotesize $1$}
      % dashed line at r
      \lpatt (0.1 0.15)
      \move (4 1) \lvec (4 -0.1)
      \htext (4 -0.5) {\footnotesize $r$}
      \lpatt ()

      % top side
      \linewd 0.07
      \move (0 1) \lvec (4 1)
    \end{texdraw}
  \end{center}
  
  A weak unit loop is by definition a null-homotopic loop together 
  with a specified null-homotopy.  To be explicit, it is a continuous
  map $\omega: T_r \to X$, where $T_r$ is the triangle
  \begin{center}
    \begin{texdraw}
      \setunitscale 20
      \linewd 0.02
      
      % coordinate axes
      \move (-0.5 0) \lvec (6 0)
      \move (0 -0.5) \lvec (0 2.5)
      
      % the filled triangle
      \move (0 0) \lvec (4 1) \lvec (0 1) \ifill f:0.86

      % top side
      \move (0 1) \lvec (4 1)
      
      \htext (-0.5 1) {\footnotesize $1$}
      \lpatt (0.1 0.15)
      \move (4 1) \lvec (4 -0.1)
      \htext (4 -0.5) {\footnotesize $r$}
      \lpatt ()
      
      \linewd 0.1
      \move  (0 1) \lvec (0 0) \lvec (4 1)
    \end{texdraw}
  \end{center}
  such that the fat sides of the triangle are mapped to the same point $x$;
  then the top side describes the null-homotopic loop.
  
  The map $U \to A$ simply sends $\omega$ to its restriction to the top side
  of the triangle.
  
  The composition of two such triangle maps, $\omega_1: T_{r_1} \to X$
  followed by $\omega_2 : T_{r_2} \to X$, is defined by skewing the second
  triangle by the linear transformation
  $$
  \fat v \mapsto \big(\begin{smallmatrix}1 & r_1 \\ 0 
  &1\end{smallmatrix}\big) \fat v$$
  and placing it next to the first triangle like this:
  \begin{center}
    \begin{texdraw}
      \setunitscale 20
      \linewd 0.02
      
      % coordinate axes
      \move (-0.5 0) \lvec (7 0)
      \move (0 -0.5) \lvec (0 2.5)
      
      % the filled triangle
      \move (0 0) \lvec (6 1) \lvec (0 1) \ifill f:0.86

      % top side
      \move (0 1) \lvec (6 1)
      
      \htext (-0.5 1) {\footnotesize $1$}
      % dashed line at r
      \lpatt (0.1 0.15)
      \move (2 1) \lvec (2 -0.1)
      \htext (2 -0.5) {\footnotesize $r_1$}
      \move (6 1) \lvec (6 -0.1)
      \htext (6 -0.5) {\footnotesize $r_2$}
      \lpatt ()
      
      \linewd 0.1
      \move  (0 1) \lvec (0 0) \lvec (2 1)
      \move ( 0 0) \lvec ( 6 1)
    \end{texdraw}
  \end{center}
  In other words, the composite is the map
  \begin{eqnarray*}
    \omega_1 \tensor \omega_2: T_{r_1+r_2} & \longrightarrow & X  \\
    \fat v & \longmapsto & 
    \begin{cases} 
      \omega_1(\fat v) & \text{ for } \fat v\in T_{r_1} \\ 
      \omega_2(\big(\begin{smallmatrix}1 & -r_1 \\ 0 
      &1\end{smallmatrix}\big) \fat v) 
      & \text{ for } \fat v\not\in T_{r_1} .
      \end{cases}
  \end{eqnarray*}
  It is clear that this composition law is strictly associative --- 
  if you wish it is because the semigroup of matrices 
  $\big(\begin{smallmatrix}1 & r \\ 0 &1\end{smallmatrix}\big)_{r>0}$
  is isomorphic to $\R_+$.  For the same reason,
  $U \to A$ is clearly compatible with composition, i.e.~is a semi-functor.
\end{blanko}
  
\newcommand{\constantx}{\ulcorner x \urcorner}

\begin{blanko}{For each $x\in X$, $U(x)$ is contractible.}
  Consider the two continuous maps
  \begin{diagram}[w=6ex]
  U(x) & \pile{\rTo^{\ell} \\ \lTo_{\constantx}} & \R_+
  \end{diagram}
  where $\ell: U(x) \to \R_+$ is the length function (length of the 
  top side of the domain triangle), and $\constantx$ is 
  the section that sends $r\in \R_+$ to the constant map on $T_r$ with value $x$.
  We claim that $\constantx \circ \ell$ is homotopic to the 
  identity map on $U(x)$, and hence $U(x)$ is contractible.  Indeed, 
  for each map $\omega : T_r \to X$, and for each $t\in [0,1]$,
  let $\omega_t$ denote the composite given by precomposing $\omega$ with 
  scaling by a factor $t$:
  \begin{center}
    \begin{texdraw}
      \setunitscale 1.3
      \move (0 35)

      \htext (6 27)  {$T_r$}
      \htext (23.5 27)  {$\longrightarrow$}
      \htext (41 27)  {$T_r$}
      \htext (59 27)  {$\longrightarrow$}
      \htext (76 27)  {$X$}
      \move (0 0) 
      \bsegment
      \lvec (0 16) \lvec (16 16) \lvec (0 0) 
      \esegment
      \move (35 0)
      \bsegment \move (0 0) \lvec (0 16) \lvec (16 16) \lvec (0 0) 
      \move (0 10) \lvec (10 10) \esegment
      
      \linewd 0.01
      \arrowheadsize l:2.5 w:1.8  
      \arrowheadtype t:V
      \move (5 11) \avec (39 6.5)
%       \move (16 16) \lvec (42 10)
%       \move (0 0) \lvec (32 0)
%       \move (0 16) \lvec (32 10)
      \htext (6 -11)  {$\fat v$}
      \htext (23.5 -10)  {$\longmapsto$}
      \htext (41 -10)  {$t\fat v$}
      \htext (59 -10)  {$\longmapsto$}
      \htext (85 -10.5)  {$\omega(t\fat v)$}
    \end{texdraw}
    \end{center}
  Now the required homotopy from $\constantx \circ \ell$ to the 
  identity of $U(x)$ is given by
  \begin{eqnarray*}
    [0,1] \times U(x) & \longrightarrow & U(x)  \\
    (t,\omega) & \longmapsto & \omega_t   .
  \end{eqnarray*}
\end{blanko}

\begin{blanko}{The maps $U \times_O A \topile A$ are homotopy equivalences.}
  The projection map is an equivalence because $U(x)$ is contractible 
  for each $x\in O$.  We shall construct a homotopy between the 
  projection map and the composition map, and hence the composition 
  map is also an equivalence.
  
  Each Moore path $\gamma:[0,r_2] \to X$ has a unique extension $\ov 
  \gamma$ to the 
  rectangle
  \begin{center}
    \begin{texdraw}
      \setunitscale 20
      \linewd 0.02
      
      % filled rectangle
      \move (0 0) \lvec (4 0) \lvec (4 1) \lvec (0 1) \ifill f:0.86
      % coordinate axes
      \move (-0.5 0) \lvec (6 0)
      \move (0 -0.5) \lvec (0 2.5)

      \htext (-0.5 1) {\footnotesize $1$}
      % r line
      \move (4 1) \lvec (4 -0.1)
      \htext (4 -0.5) {\footnotesize $r_2$}

      % top side
      \linewd 0.06
      \move (0 1) \lvec (4 1)
    \end{texdraw}
  \end{center}
  constant in the vertical direction. 
  Given a Moore loop $\omega: T_{r_1} \to X$
  based at $\gamma(0)$, define the composite 
  of $\omega$ with $\ov\gamma$ to be the map $\omega\boxtimes\ov\gamma$
  with domain
  \begin{center}
    \begin{texdraw}
      \setunitscale 20
      \linewd 0.02
      
      % the filled triangle
      \move (0 0) \lvec (4 0) \lvec (6 1) \lvec (0 1) \ifill f:0.86
      
      % coordinate axes
      \move (-0.5 0) \lvec (7 0)
      \move (0 -0.5) \lvec (0 2.5)
      
      % top side
      \move (0 1) \lvec (6 1) \lvec (4 0)
      
      \htext (-0.5 1) {\footnotesize $1$}
      % dashed line at r
      \lpatt (0.1 0.15)
      \move (2 1) \lvec (2 -0.1)
      \htext (2 -0.5) {\footnotesize $r_1$}
      \htext (4 -0.5) {\footnotesize $r_2$}
      \move (6 1) \lvec (6 -0.1)
      \htext (6 -0.5) {\footnotesize $r_1\!+\!r_2$}
      \lpatt ()
      \move (4 0.05) \lvec (4 -0.07)
      
      \linewd 0.1
      \move  (0 1) \lvec (0 0) \lvec (2 1)
%       \move ( 0 0) \lvec ( 6 1)
    \end{texdraw}
  \end{center}
given by
\begin{eqnarray*}
    \fat v & \longmapsto & 
    \begin{cases} 
      \omega(\fat v) & \text{ for } \fat v\in T_{r_1} \\ 
      \ov\gamma(\big(\begin{smallmatrix}1 & -r_1 \\ 0 
      &1\end{smallmatrix}\big) \fat v) 
      & \text{ for } \fat v\not\in T_{r_1}
      \end{cases}
  \end{eqnarray*}
  This defines a continuous map
  $$U \times_O A \to B
  $$
  where $B$ is the space of maps to $X$ from such trapezia.
  Now the composition map $U \times_O A \to A$ is given
  by postcomposing this map with `restriction to the $(y\!=\! 1)$-line'
  \begin{eqnarray*}
    U \times_O A & \longrightarrow & \phantom{xi}B\phantom{wi} \ 
    \longrightarrow \ \ \ \ A  \\
    (\omega,\gamma) & \longmapsto & \omega\boxtimes\ov\gamma  \ \longmapsto \ 
    (\omega\boxtimes\ov\gamma) \!\mid_{y=1} \;=\; \omega\tensor \gamma .
  \end{eqnarray*}
  The required homotopy from projection ($t=0$)
  to composition ($t=1$) is now given by
  \begin{eqnarray*}
    [0,1] \times U \times_O A & \longrightarrow & \phantom{xi}B\phantom{wi} \  
    \longrightarrow \ \ \ \ A  \\
    (t,\; \omega,\gamma) & \longmapsto & \omega\boxtimes\ov\gamma  \ \longmapsto \ 
    (\omega\boxtimes\ov\gamma) \!\mid_{y=t}
  \end{eqnarray*}
  consisting in sliding up the support line $y=t$:
  \begin{center}
    \begin{texdraw}
      \setunitscale 20
      \linewd 0.02
      
      % the filled triangle
      \move (0 0) \lvec (4 0) \lvec (6 1) \lvec (0 1) \ifill f:0.86
      
      % coordinate axes
      \move (-0.5 0) \lvec (7 0)
      \move (0 -0.5) \lvec (0 2.5)
      
      % top side
      \move (0 1) \lvec (6 1) \lvec (4 0)
      
%       \htext (-0.5 1) {\footnotesize $1$}
      % dashed line at r
      \lpatt (0.1 0.15)
      \move (2 1) \lvec (2 -0.1)
      \htext (2 -0.5) {\footnotesize $r_1$}
      \htext (4 -0.5) {\footnotesize $r_2$}
      \move (6 1) \lvec (6 -0.1)
      \htext (6 -0.5) {\footnotesize $r_1\!+\!r_2$}
      \lpatt ()
      \move (4 0.05) \lvec (4 -0.07)

      \move  (0 1) \lvec (0 0) \lvec (2 1)
%       \linewd 0.1
      \move (-0.2 0.55) \lvec (6.5 0.55)
      \htext (-0.5 0.55) {\footnotesize $t$}
    \end{texdraw}
  \end{center}
  To be pedantic with the domain, we should say that the path 
  $(\omega\boxtimes\ov\gamma)\!\mid_{y=t}$ is defined to be the map
  \begin{eqnarray*}
    [0,tr_1+r_2] & \longrightarrow & X  \\
    s & \longmapsto & (\omega\boxtimes\ov\gamma)(\begin{smallmatrix}s \\ t 
  \end{smallmatrix})
  \end{eqnarray*}
%   
% 
%   
%   The composition map is defined as
%   \begin{eqnarray*}
%     U \times_O A & \longrightarrow & A  \\
%     (G,\gamma) & \longmapsto & G \!\mid_{(y=1)} \;\#\; \gamma .
%   \end{eqnarray*}
%   The homotopy from the projection map can be written sleekly as
%   \begin{eqnarray*}
%     [0,1] \times U \times_O A & \longrightarrow & A  \\
%     (t, \,G,\gamma) & \longmapsto & G \!\mid_{(y=t)} \;\#\; \gamma ,
%   \end{eqnarray*}
%   but this may be seen as cheating since $G\!\mid_{t=0}$ is not a 
%   valid Moore path.

  This finishes the proof of the proposition.
  \qed
\end{blanko}

%%%%%%%%%%%%%%%%%%%%%%%%%%%%%%%%%%%%%%%%%%%%%%%%%%
\subsection*{Monoidal model categories}
%%%%%%%%%%%%%%%%%%%%%%%%%%%%%%%%%%%%%%%%%%%%%%%%%%

The next example is of a {\em fair monoidal coloured category}, i.e.,
a fair monoid in $\CCat$ (cf.~\ref{CCat}).  To make sense of this it
must be specified in which sense the category $\CCat$ is a coloured
category.  There are several possible colour structures --- the
crucial desired property is that the Dwyer-Kan simplicial localisation
functor $L: \CCat \to \sCat$ \cite{Dwyer-Kan:simplicial-localization}
should be colour-preserving.  For
simplicity we take this as the definition: a (colour-preserving)
functor $F : (C,W) \to (C',W')$ between coloured categories is called
an {\em equifunctor} if $LF$ is an equivalence of simplicial
categories.  This means that $\pi_0 LF : \pi_0 LC \to \pi_0 LC'$ is an
equivalence of categories and for each pair of objects $x,y\in C$
the map $LC(x,y) \to LC'(Fx, Fy)$ is a weak equivalence of simplicial
sets.

\begin{blanko}{Cofibrant objects of a monoidal model category.}
  Let $(M,\tensor, I)$ be a monoidal model category in the sense of
  Hovey~\cite{Hovey:model}.  This means that $M$ is at the same time a
  model category and a monoidal category and that the following
  compatibility conditions are met:
  \begin{punkt-i}
    \item   $\tensor$ preserves weak equivalences between cofibrant objects; 
  
    \item for every weak equivalence $Z \isopil I$ with $Z$ cofibrant,
    and for every cofibrant object $X$, the composite $Z\tensor X \to
    I \tensor X \to X$ is a weak equivalence.
  \end{punkt-i}
  Now if $I$ happens to be cofibrant then the subcategory of cofibrant
  objects $(M^c, \tensor , I)$ is a genuine monoidal coloured category
  with unit (i.e., a monoid in $\CCat$), and induces a monoidal
  structure on the homotopy category $ \Ho(M)$.  If $I$ is not
  cofibrant then $(M^c,\tensor)$ will not have a unit object.
  Instead:
  
  \begin{lemma}
  The category of cofibrant objects $M^c$ is naturally a fair monoidal
  coloured category (i.e., a fair monoid in $\CCat$). 
  The weak units are the equivalences $Z \isopil I$ with $Z$ cofibrant.
  \end{lemma}
  
  The functor $\grosop\to\CCat$ is given by taking $O$ to be
  singleton, $A$ to be the coloured category $M^c$, and $U$ the coloured category
  whose objects are weak equivalences $Z \isopil I$ with $Z$ cofibrant, and
  whose arrows are triangles $Z'\isopil Z \isopil I$.  The multiplication
  functor $\tensor$ clearly turns $A$ and $U$ into semi-monoids as required,
  hence the Segal condition is satisfied.  To check the weak identity arrow
  axiom, it must first be shown that $U$ is contractible; this is done by
  comparing it with the category $M^c/I'$ of cofibrant objects over a fixed
  cofibrant replacement $I' \isopil I$.  This category has a terminal
  object and is therefore contractible.  Second, it must be shown that the
  two functors $U \times A \pile{\rTo^m \\
  \rTo_p} A$ are equifunctors.  The projection $p$ is an equifunctor
  because $U$ is contractible, and $m$ is an equifunctor because the equivalences
  of (ii) assemble into a natural transformation $m\Rightarrow p$.
\end{blanko}
  
%   \begin{blanko}{Examples}
%     of monoidal model categories whose unit is not cofibrant include
%     categories of $S$-modules \cite{Schwede-Shipley:9801}\ldots
%     
%   \end{blanko}
  
\begin{blanko}{Spitzweck's monoidal model categories with pseudo-unit.}
  The fair monoid structure on $M^c$ does not depend on the fact that $I$
  is a genuine unit: it is enough to have maps $I \tensor X \to X$ such
  that for every $Z\isopil I$ with $Z$ cofibrant, the composite $Z\tensor X
  \to I \tensor X \to X$ is a weak equivalence, for cofibrant $X$.  Such
  model categories have been studied by Spitzweck~\cite{Spitzweck:0101102}
  under the name `monoidal model categories with pseudo-unit'.  The lemma
  holds also in this situation.

  Such structures arise in connection with modules and algebras over
  $E_\infty$-operads.  Let $\mathscr{L}$ denote the linear isometries
  operad --- it is an unital $E_\infty$-operad (see
  \cite{Elmendorf-Kriz-Mandell-May} for more information).  Consider the
  category $S\kat{-Mod}$ of modules over the differential graded algebra $S
  \df \mathscr{L}(1)$.  It was shown in
  K\v{r}\'{\i}\v{z}-May~\cite{Kriz-May} that $S\kat{-Mod}$ can be equipped
  with an associative tensor product (analogous to the one constructed in
  \cite{Elmendorf-Kriz-Mandell-May} for spectra), $X\boxtimes Y \df
  \mathscr{L}(2) \tensor_{S\tensor S} X \tensor Y$, but there is no unit
  for $\boxtimes$.  A similar construction works for modules over a given
  $S$-algebra.
  
  Now if $M$ is a model category satisfying some mild technical conditions,
  all these notions and constructions make sense in $M$, and
  Spitzweck~\cite{Spitzweck:0101102} shows that the category of $S$-modules
  in $M$ is a monoidal model category with pseudo-unit, and that if $A$ is
  a cofibrant $S$-algebra the same is true for the category of $A$-modules
  in $M$.  Hence, in each case the category of cofibrant objects becomes a
  fair monoidal coloured category.
\end{blanko}

\subsection*{Cobordism categories}
%%%%%%%%%%%%%%%%%%%%%%%%%%%%%%%%%%%%%%%%%%%%%%%%%%
\label{sub:Cob}

The final example in this exposition is an example of a  fair Tamsamani
$2$-category.  That is, a $\grosop$-diagram in $\Cat$ whose 
Segal maps are only equimorphisms (and satisfying the other axioms for 
a fair category).

\begin{blanko}{Classical cobordism categories.}
  (See~\cite{Kock:FA-2DTQFT} for all details.)
  Classically, the category $\nCob$ of oriented $n$-dimens\-ional
  cobordism is the category whose objects are closed oriented
  $(n-1)$-manifolds, and whose morphisms from $\Sigma_0$ to $\Sigma_1$
  are equivalence classes of oriented cobordisms whose `in-boundary' 
  is $\Sigma_0$
  and whose `out-boundary' is $\Sigma_1$.  Here two cobordisms are 
  considered equivalent if there is a diffeomorphism between them that
  induces the identity map on the boundaries.  It is necessary to pass 
  to this quotient because given two manifolds with matching 
  boundaries the gluing is not canonical:  the gluing does exist and
  all possible gluings are diffeomorphic, but there is no universal 
  property and no unique comparison diffeomorphism.  Note also that
  identity cobordisms cannot exist before dividing out by 
  diffeomorphisms, since the only true identities would be the cylinders 
  of height zero, and they are not $n$-manifolds.
\end{blanko}

\begin{blanko}{Non-algebraic composition law.}
  Instead of choosing a specific composition for two given cobordisms,
  which as explained would be an artificial choice anyway, we can
  indicate the space of all possible compositions.  In other words, we
  will define a functor $X:\Dmono\op \to \Cat$ and argue that it
  satisfies the weak Segal condition.  The set of objects $X_0$ is the
  set of all closed oriented $(n-1)$-manifolds.  The category $X_1$
  has as objects the oriented cobordisms; the two maps $X_1 \topile
  X_0$ associate to a given cobordism its in-boundary and its
  out-boundary.  The arrows of the category $X_1$ are the homotopy classes
  of diffeomorphisms that restrict to the identity on the boundaries.
  To be more precise, a cobordism from $\Sigma_0$ to $\Sigma_1$ is a
  quintuple $(M;\Sigma_0, \Sigma_1; \sigma_0,\sigma_1)$ where
  $\sigma_0$ is a diffeomorphism from $\Sigma_0$ onto the in-boundary
  of the $n$-manifold $M$, and $\sigma_1$ is a diffeomorphism from
  $\Sigma_1$ onto its out-boundary.  The arrows in the category
  between $(M;\Sigma_0, \Sigma_1; \sigma_0,\sigma_1)$ and
  $(M';\Sigma_0, \Sigma_1; \sigma_0', \sigma_1')$ are then the
  homotopy classes of diffeomorphisms $M\isopil M'$ compatible with the
  $\sigma_i$.
  
  So far the description is classical.  Now instead of having a
  composition map, we specify the category $X_2$: its objects are
  cobordisms with a specified decomposition.  Precisely they are
  septuples $(M; \Sigma_0,\Sigma_1,\Sigma_2;
  \sigma_0,\sigma_1,\sigma_2)$, where $\sigma_0$ is a diffeomorphism
  from $\Sigma_0$ onto the in-boundary of $M$; $\sigma_2$ is a
  diffeomorphism from $\Sigma_2$ onto the out-boundary of $M$, and
  $\sigma_1$ is an embedding of $\Sigma_1$ into $M$ splitting it into
  two cobordisms, one $M'$ from $\Sigma_0$ to $\Sigma_1$, and another
  $M''$ from $\Sigma_1$ to $\Sigma_2$.  A morphism between two such
  objects is just a homotopy class of diffeomorphisms of the cobordism
  compatible with the sigmas.
  
  Now we have to specify the three face maps $X_2 \pile{\rTo\\ \rTo \\
  \rTo} X_1$: these are: return $M''$, return $M'$, or return the whole
  $M$ forgetting the splitting data $(\Sigma_1,\sigma_1)$.  In
  general, the category $X_k$ is defined to be the category of tuplets
  $(M;\Sigma_i,\sigma_i)_{0\leq i \leq k}$ and diffeomorphisms preserving the
  submanifolds as explained for $k=2$.  It is clear that all these
  categories form a functor $\Dmono\op\to\Cat$.  To verify that the
  Segal holds condition we must show that such a subdivided cobordism
  is determined up to diffeomorphism by the pieces it is made up of.
  This statement follows from the fact that it is possible to glue up
  to diffeomorphism.
\end{blanko}

\begin{blanko}{Cylinders as weak units.}
   There is no way of extending this functor $X:\Dmono\op\to\Cat$
  to $\Delta$ to get a simplicial category.  But there is a natural
  extension to $\gros$.  To begin with, send \intextBbB\ to the category of
  straight cylinders: the objects are the cobordisms $\Sigma\times I$
  (and more generally, cobordisms $\Sigma\to\Sigma$ equipped with a
  diffeomorphism to such a straight cylinder), and the arrows are
  homotopy classes of boundary preserving diffeomorphisms induced from
  diffeomorphisms of intervals. Clearly this category is contractible.
  There is a functor from this category to $X_1$ consisting in 
  forgetting the straight structure. 
  
  More generally, an object $K\in\gros$ is sent to the category of all 
  split cobordisms such that the pieces corresponding to links in 
  $K$ are straight cylinders.  Again the vertical maps are seen to be 
  equivalences.  Altogether we have defined a colour-preserving 
  functor
  $$
  X : \grosop\to \Cat
  $$
  with $X\!_{\intextB}$ discrete, and satisfying the weak Segal condition, and 
  hence a fair Tamsamani category in $\Cat$.  In fact the hom cats are all 
  groupoids.
\end{blanko}

\begin{prop}\label{CobProp}
  Oriented $n$-cobordisms naturally assemble into a fair Tamsamani
  $2$-category, for which the straight cylinders are weak identity
  arrows.
\end{prop}

\begin{BM}
  The above construction is clearly $2$-truncated.  It is possible to
  avoid this truncation: instead of getting a fair Tamsamani category in
  groupoids, the result is a fair Tamsamani category in simplicial
  sets, i.e.~a fair Segal category, reflecting one level of cobordism
  structure and reflecting homotopy for all higher levels.  Simplicial
  representations of such fair Segal categories of cobordisms might be
  an interesting alternative approach to the idea of extended
  topological quantum field theories advocated by
  Lawrence~\cite{Lawrence:extended} and
  Baez-Dolan~\cite{Baez-Dolan:TQFT-higher}.
\end{BM}

% 
% \bigskip
% 
% 
% Clearly the examples above have a lot in common: in each case
% there is actually a formal unit somewhere, but since it is not of the good 
% type it is replaced by a contractible space of suitable approximations;
% in each case these approximation objects come equipped with a comparison to 
% the `excluded' strict unit, and in particular, the space $U$ is not a 
% subspace of the space of arrows $A$.

\renewcommand{\thesection}{\Alph{section}}
\setcounter{section}{1}
\setcounter{lemma}{0}
\addcontentsline{toc}{section}{\numberline{}Appendix: Discrete objects}

%%%%%%%%%%%%%%%%%%%%%%%%%%%%%%%%%%%%%%%%%%%%%%%%%%
\section*{Appendix: Discrete objects}
%%%%%%%%%%%%%%%%%%%%%%%%%%%%%%%%%%%%%%%%%%%%%%%%%%
\label{Sec:discrete}

In the main text, three crucial notions are `discrete objects', 
`fibre products over discrete objects', and `equimorphisms', and 
the important functors are those that preserve these notions.
For flexibility, `discrete objects' is considered a structure to be
specified, not a property;  the main reason is that we want to say
that \intextB\ is the only discrete object in $\gros$ (or in $\Delta$),
whereas this case is not covered by the notion of standard discrete 
objects described below.

Just for the notions of $\SS$-category and
fair $\SS$-category to make sense, the only requirement on the
discrete objects in $\SS$ is that $\SS$ should admit fibre products 
over discrete objects.  However, in order to get a 
reasonable theory, and in particular to get a good notion of
equimorphisms in the categories of $\SS$-categories and fair 
$\SS$-categories, certain features of the discrete objects are needed,
requiring in turn certain closure properties of $\SS$, and finally
the notion of equimorphism in $\SS$ must be compatible with these notions.

% 
% \bigskip
% 
% There are three aspects involved:
% 
% --- closure properties
% 
% --- discrete objects
% 
% --- colours
% 
% and between each of these three aspects some obvious compatibilities 
% should hold.  We treat the three aspects in that order:

\begin{blanko}{Closure properties.}\label{closure}
  We require of $\SS$ that it has all sums, and these
  should be
  disjoint and universal (cf.~SGA 4.1 \cite{SGA4.1}, Exp.~II,
  Def.~4.5).  Recall that a sum $A=\coprod_{i\in I} A_i$ is
  called disjoint when each structure map $A_i \to A$ is a
  monomorphism and admits all pullbacks, and for each $i\neq j$ the
  pullback $A_i \times_A A_j$ is an initial object in $\SS$.
  Universal means that the pullback of a sum diagram is again a
  sum diagram.
  
  Next we require $\SS$ to have all finite products.  It follows 
  automatically from the sum requirements
  that finite products distribute over sums,
  $$
  A \times \smallcoprod{i\in I}{} B_i \simeq \coprod_{i\in I} (A 
  \times B_i) .
  $$
  % PROOF: Set $S= A\times \smallcoprod{i\in I}{} B_i$.  For each
  % $i\in I$ we have cartesian squares
%   \begin{diagram}[w=7ex,h=4.5ex,tight]
%   S_i \SEpbk& \rTo  & B_i  \\
%   \dTo  &    & \dTo  \\
%   S \SEpbk & \rTo  & \smallcoprod{i\in I}{} B_i \\
%   \dTo && \dTo \\
%   A & \rTo & *
%   \end{diagram}
%   where $S_i$ is the pullback by definition.  By universality, $S = 
%   \smallcoprod{i\in I}{} S_i$, but on the other hand the outer square
%   is also a pullback, so we have $S_i = A \times B_i$.
\end{blanko}

\begin{blanko}{Standard discrete objects.}\label{SDO}
  For a category $\SS$ with the closure properties of \ref{closure},
  the {\em standard discrete objects structure} is the case where the
  discrete objects are the sums of terminal objects $*$, and
  the following two conditions are satisfied:
    
    [DO1] \ The discrete objects 
    form a full reflective subcategory. That is, the
    discrete-objects functor 
    \begin{eqnarray*}
	\delta : \Set & \longrightarrow &  \SS  \\
	I & \longmapsto & \coprod_{i\in I} *
    \end{eqnarray*}
    is fully faithful and has a left adjoint denoted $\pi_0$ (the components 
    functor).
    (Note that there is 
    always a right adjoint: $A \mapsto \Hom_\SS (*, A)$.)

    [DO2]  \ The components functor $\pi_0$ preserves finite products.
\end{blanko}

The second axiom simply expresses full compatibility of the adjunction 
with the stipulated sums and products --- indeed it is already
automatic that $\pi_0$ preserves sums and that $\delta$ preserves
products, and $\delta$ also preserves sums by construction.

\begin{BM}
  The description of the discrete-objects functor $\delta: \Set \to
  \SS$ should not be taken too literally as a specific functor defined
  in terms of an arbitrary choice of $\coprod$ and $*$ in $\SS$;
  rather it is characterised by a universal property (preservation of
  sums and finite products), and any object isomorphic to a 
  discrete object should be discrete again.  
  In fact the category $\Set$ could be regarded as a mere place holder
  for a genuine subcategory-of-discrete-objects-in-$\SS$, i.e., 
  derived from $\SS$.
\end{BM}

% For example, if literally we say that a discrete category
% is one where the only arrows are the identity arrows and that a discrete 
% simplicial set is a constant presheaf on $\Delta$, then the nerve
% of a discrete category would not be discrete, since the set of objects
% of a category is not identical to its set of identity arrows.

\begin{blanko}{Decomposition.}\label{decomp}
  Given an arrow $A \to \delta I$ to a discrete object; for each $i\in
  I$ define $A_i$ to be the fibre product
    \begin{diagram}[w=5ex,h=4.5ex,tight]
    A_i\SEpbk & \rTo  & A  \\
    \dTo  &    & \dTo  \\
    *  & \rTo_i  & \delta I .
    \end{diagram}
  It follows from the universality of sums that there is
  a natural isomorphism
  $$
  A \simeq \coprod_{i\in I} A_i .
   $$
\end{blanko}

\begin{blanko}{Fibre products over discrete objects.}\label{exist-fibr}
  More generally, we get the
%   The decomposition axiom 
% %   has several implications: the initial 
% %   object $\emptyset = \delta\emptyset$ has the property that any 
% %   arrow $A \to \emptyset$ is an isomorphism.  (Indeed, decompose
% %   $A = \coprod_{x\in \emptyset} A(x)$.)  Furthermore, 
 existence of all fibre products
  over discrete objects: 
  $$
  A \times_{\delta I} B \simeq \coprod_{i\in I} A_i \times B_i.
  $$
  % Here we use naturality of decomposition: given an arrow $f : A 
  % \to B$ over a discrete object, then $f= \coprod f_i : A_i \to B_i$.
  % Now show that $\coprod_{i\in I} A_i \times B_i$ has the universal
  % property of the product: given $Z$ with arrows to $A$ and $B$
  % (coinciding down in $\delta I$), decompose over $I$ and use the
  % universal property of the cartesian product to get a unique map
  % $f_i : Z_i \to A_i \times B_i$.  The sum of all these maps
  % is what we are looking for.  This map is unique, because it must
  % decompose, and the components are unique\ldots
  % 
  It then follows from 
  axiom DO2 that
  $\pi_0$ preserves fibre products over discrete objects.
  
  Also, sums commute with fibre products over discrete 
  objects:
  \begin{equation}\label{commute}
  \coprod_{i\in I} A_i \times_{\delta X_i} B_i \ \simeq \
  \smallcoprod{i\in I}{} A_i
  \underset{\smallcoprod{i\in I}{}\delta X_i}{\times} \smallcoprod{i\in I}{}B_i.
  \end{equation}
  
  Let us see how this goes by a quick  example computation which also shows the 
  importance of the rule: the Segal condition is preserved under 
  sums. Given two fibre products
  $$ 
  A_r \simeq A_p \times_{A_0} A_q \qquad
  B_r \simeq B_p \times_{B_0} B_q
  $$
  with $A_0$ and $B_0$ discrete
  (e.g.~two Segal conditions!), put $C_i \df A_i \coprod B_i$ 
  ($i=0,p,q,r$). 
  Then the statement is that $C_r \simeq C_p 
  \times_{C_0} C_q$.  
  Indeed, in the fibre product $C_p
  \times_{C_0} C_q \simeq \big( A_p \coprod B_p\big) \times_{A_0\coprod B_0} \big( 
  A_q \coprod B_q\big)$, decompose over the discrete 
  object $A_0 \coprod B_0$ to get
  $$
  \simeq \coprod_{x\in A_0\coprod B_0} ( A_p \coprod B_p)_x \times ( 
  A_q \coprod B_q)_x .
  $$
  Now write this as the sum of two sums, and use the disjointness
  axiom to remove half of the summands
  $$
  \simeq 
  \smallcoprod{x\in A_0}{} \big((A_p)_x \times (A_q)_x\big)
  \coprod
    \smallcoprod{x\in B_0}{} \big((B_p)_x \times (B_q)_x\big)
  \simeq 
  A_r \coprod B_r  \simeq C_r.
  $$
\end{blanko}

\begin{blanko}{Examples.}
  The archetypical example is $\Top$.  The discrete spaces are
  standard discrete objects, and $\pi_0$ is the usual components
  functor.  In $\Set$, every object is discrete, and $\delta$ and
  $\pi_0$ are the identity functors.  The key examples for the present
  purposes are $\sSet$ and $\Cat$, which we analyse a little further.
\end{blanko}
  
\begin{blanko}{Simplicial sets.}\label{pi0}
  If $D$ is a small category, the presheaf category $\Cat(D\op,\Set)$
  satisfies the closure properties in \ref{closure} as well as axiom
  DO1 in \ref{SDO} (here $\delta$ is the constant-presheaf functor,
  and $\pi_0$ returns the colimit of a given $D\op$-diagram).
  However, axiom DO2 is not in general satisfied.  Satisfying axiom
  DO2 is a crucial property of a category $D$ supposed to serve as
  base category for combinatorial topology.  (This was observed by
  Grothendieck~\cite{Grothendieck:stacks} and has also been stressed
  by Lawvere.)  The category $\Delta$ has this property, which is to
  say that $\sSet$ has standard discrete objects.  To see this, note
  that if $X$ is a simplicial set then $\pi_0(X)$ is described as the
  quotient of $X_0$ by the equivalence relation defined by identifying
  two $0$-cells if they can be connected by `zigzags' of $1$-cells.
  The key point here is the existence of degeneracy maps,
  i.e.~reflexivity of the equivalence relation.  (In contrast $\Dmono$
  does not have the property: consider the $\Dmono\op$-diagram $X$ where
  $X_0$ is a two-element set, $X_1$ is singleton, $X_0\topileback X_1$
  are the two inclusions, and $X_n=\emptyset$ for $n\geq 2$.  Then
  $\pi_0(X\times X) \not\simeq \pi_0(X) \times \pi_0(X)$.  The fat delta
  $\gros$ does not have the property either (the previous example can
  be used again), but as we shall see, in the setting of coloured
  categories, it {\em does} have the property.)
\end{blanko}

\begin{blanko}{Categories.}
  The full subcategory $\Cat\subset\sSet$ is closed under sums
  and finite products, and enjoys the closure properties of
  \ref{closure}.  Indeed, the Segal condition is given in terms of
  fibre products over discrete objects, and these commute with sums
  and finite products.

  The importance of axiom DO2 is that it furthermore allows for the
  discrete objects structure to descend from $\sSet$ to
  $\Cat$.  Indeed, since $\pi_0^{\sSet}$ (and also $\delta^{\sSet}$) 
  preserves fibre products over 
  discrete objects,  the discrete-objects adjunction for $\sSet$
  restricts to an adjunction for $\Cat$:
  \begin{diagram}[w=6.5ex,h=4ex,tight]
  \sSet & \pile{\lTo^{\delta^{\sSet}} \\ \rTo_{\pi_0^{\sSet}}}  & \Set  \\
  \cup &    &  \parallel \\
  \Cat & \pile{\lTo^{\delta^{\Cat}} \\ \rTo_{\pi_0^{\Cat}}}  & \Set
  \end{diagram}
  describing the standard discrete objects structure on $\Cat$.
\end{blanko}

These two examples readily generalise to the case where $\Set$ is 
replaced by a general category $\SS$ with standard discrete objects:

\begin{lemma}\label{Cat-Delta-S}
  If $\SS$ has standard discrete objects then so has $\Cat(\Delta\op,\SS)$.
\end{lemma}

\newdiagramgrid{gridge}{1.8,1}{1}

\begin{dem}
  The standard-discrete-objects adjunction $\pi_0 
  \isleftadjointto \delta$ for $\SS$
  induces an adjunction
  \begin{diagram}[w=12ex,h=4.5ex,tight]
  \Cat(\Delta\op,\SS) & \pile{\lTo^{\delta\lowerstar} 
  \\ \rTo_{\pi_0{}\lowerstar}} & \Cat(\Delta\op,\Set)=\sSet
  \end{diagram}
  by postcomposition with $\delta$ and $\pi_0$.
  Since $\delta$ and $\pi_0$ preserve sums and finite products,
  and since sums and finite products are computed 
  point-wise, the induced functors $\delta\lowerstar $ and 
  $\pi_0{}\lowerstar $ again preserve sums and finite 
  products.  The composite adjunction
  \begin{diagram}[w=5ex,h=3.5ex,grid=gridge,tight]
  \Cat(\Delta\op,\SS) & \pile{\lTo 
  \\ \rTo} &\sSet& \pile{ \lTo \\ \rTo} &\Set
  \end{diagram}
  describes the standard discrete objects structure on
  $\Cat(\Delta\op,\SS)$ (the discrete objects are the constant
  presheaves with discrete value).
% Further details: the claim that $\pi_0{}\lowerstar$ preserves 
% products: let $F$ and $G$ be simplpicial objects in $\SS$.  Their 
% product is the simplicial object whose $k-cells are $(F\times G)_k =
% F_k \times G_k$.  Now postcompose with $\pi_0$.  Clearly
% $(\pi_0 F \times \pi_0 G)_k = \pi_0( F \times G)_k$\ldots
\end{dem}
  \begin{lemma}\label{discSCat}
   If $\SS$ has standard discrete objects, then $\SCat$
   has standard discrete objects.
\end{lemma}

\begin{dem}
  First notice that the full subcategory $\SCat\subset
  \Cat(\Delta\op,\SS)$ inherits sums and \mbox{finite} products from
  the ambient category, and hence satisfies the closure properties of
  \ref{closure}.  Now the two adjoint functors $\Cat(\Delta\op,\SS)
  \pile{\lTo \\ \rTo} \Set $ of the previous lemma \mbox{preserve}
  discrete objects and fibre products over discrete objects.  Hence
  they restrict to two functors $\SCat \pile{\lTo \\ \rTo}\Set$, and
  since $\SCat\subset \Cat(\Delta\op,\SS)$ is full this pair of
  functors forms again an adjunction, and since sums and finite
  products of categories are computed as simplicial sets, these two
  functors also preserve sums and finite products.
\end{dem}

%%%%%%%%%%%%%%%%%%%%%%%%%%%%%%%%%%%%%%%%%%%%%%%%%%
\subsection*{Discrete objects and colours}
%%%%%%%%%%%%%%%%%%%%%%%%%%%%%%%%%%%%%%%%%%%%%%%%%%

\begin{blanko}{Colours and standard discrete objects.}\label{CSDO}
  Assume $\SS$ has standard discrete objects.  We shall now define
  what it means for a colour structure on $\SS$ to be compatible with
  the standard discrete objects.  There are compatibility conditions
  with respect to the sums and finite products, and a compatibility
  condition with respect to the discrete-objects adjunction:

  --- Stability under sums and finite products: If $(f_i)_{i\in
  I}$ is a family of equimorphisms then $\coprod_{i\in I} f_i$ is
  again an equimorphism, and given equimorphisms $f_1,\ldots,f_n$,
  then the product $\prod_{i=1}^n f_i$ is again an equimorphism.

  --- Preservation under $\pi_0$: if $f$ is an equimorphism in $\SS$
  then $\pi_0(f)$ is an equimorphism in $\Set$, i.e.~a bijection.
  Note that it is automatic that $\delta$ preserves equimorphisms.
  The condition on $\pi_0$ implies that the restriction of $\delta$ to
  the subcategory $\kat{Bij}$ of sets and bijections is again a full
  embedding $\kat{Bij} \to W$.
\end{blanko}

\begin{blanko}{Examples.}
  The categories $\sSet$, $\Top$, $\Set$, and $\Cat$, with the usual
  notions of discrete objets and equimorphisms, all fit into the
  framework above.  Here $\Set$ has bijections as equimorphisms, and
  it follows that in each case the functors $\pi_0$ and $\delta$
  themselves preserve sums, finite products, and equimorphisms.  The
  geometric realisation functor $\sSet\to\Top$ as well as the nerve
  functor $\Cat\to\sSet$ are also examples of functors preserving
  sums, finite products, and equimorphisms.
\end{blanko}

The components functor $\pi_0: \Cat(\grosop,\Set)\to\Set$ does not
preserve finite products, for the same reason as mentioned for
$\Dmono$ in \ref{pi0}.  However, $\gros$ was not designed with
presheaves in mind, but rather {\em restricted presheaves} with
respect to preservation of colours: we are interested in the full
subcategories $\CCat(\grosop,\Set)\subset \Cat(\grosop,\Set)$ and
$\CCat(\grosop,\SS)\subset \Cat(\grosop,\SS)$, and in here the full
subcategories of fair $\SS$-categories and fair $\Set$-categories.
   
\begin{lemma}
  If $\SS$ has standard discrete objects with compatible colours,
  then $\CCat(\grosop,\SS)$ has standard discrete objects.
\end{lemma}

\begin{dem}
  Since the equimorphisms in $\SS$ are assumed to be stable under
  sums and finite products, it follows that the full subcategory
  $\CCat(\grosop,\SS) \subset \Cat(\grosop,\SS)$ satisfies the
  closure properties of \ref{closure}. 

  Just like in \ref{Cat-Delta-S}, 
  the standard-discrete-objects adjunction $\pi_0 
  \isleftadjointto \delta$ for $\SS$
  induces an adjunction
  \begin{diagram}[w=12ex,h=4.5ex,tight]
  \CCat(\grosop,\SS) & \pile{\lTo^{\delta\lowerstar} 
  \\ \rTo_{\pi_0{}\lowerstar}} & \CCat(\grosop,\Set)
  \end{diagram}
  by postcomposition with $\delta$ and $\pi_0$.  This time, in order
  for this to work it is crucial that $\delta$ and $\pi_0$ preserve
  equimorphisms.  Since $\delta$ and $\pi_0$ also preserve sums and finite
  products, the induced functors $\delta\lowerstar $
  and $\pi_0{}\lowerstar $ again preserve sums and finite products.
  
  In Proposition~\ref{sSet=FairSSet} we established an adjoint
  equivalence $\CCat(\grosop,\Set) \pile{\lTo \\ \rTo} \sSet$, and it
  is easy to check that each of these two adjoint functors preserves
  sums and finite products.  Now the standard discrete objects
  in $\CCat(\grosop,\SS)$ are described by the composite adjunction
  $$
  \CCat(\grosop,\SS)
   \pile{\lTo \\ \rTo}
  \CCat(\grosop,\Set)
  \pile{\lTo \\ \rTo} \sSet \pile{\lTo \\ \rTo} \Set .
  $$
\end{dem}
% Note that at this point it makes no sense to ask if these functors
% are colour preserving, because we have not defined any colour 
% structure on $\CCat(\grosop,\SS)$.  Rather, if you choose a colour
% structure on $\CCat(\grosop,\SS)$ then you should choose it 
% compatible with the standard discrete objects structure described in 
% the lemma.

\begin{prop}\label{SFairCat-SDO}
   If $\SS$ has standard discrete objects with compatible colours, 
   then $\kat{S-FairCat}$
   has standard discrete objects.
\end{prop}

\begin{dem}
  Just like in \ref{discSCat}, the full subcategory
  $\kat{S-FairCat}\subset \CCat(\grosop,\SS)\subset \Cat(\grosop,\SS)$ 
  inherits sums and finite products from the ambient category,
  and hence satisfies the closure properties of \ref{closure}.
  Now the adjoint functors of the previous lemma all preserve
  discrete objects, 
  sums and finite products, and hence
   fibre products over discrete objects, so they 
  restrict to adjoint functors
%   \begin{diagram}[w=10ex,h=4.5ex,tight]
%   \CCat(\grosop,\SS)&\pile{\lTo \\ \rTo}& \CCat(\grosop,\Set)  &
%  \pile{\lTo \\ \rTo}& \sSet& \pile{\lTo \\ \rTo}& \Set \\
%    \cup &    & \cup &&\cup&&\parallel  \\
   $$
  \kat{S-FairCat}  \pile{\lTo \\ \rTo} \kat{Set-FairCat} 
  \pile{\lTo \\ \rTo} \Cat \pile{\lTo \\ \rTo} \Set
  $$
%   \end{diagram}
  (which again preserve sums and finite products).
  This describes the standard
  discrete objects structure on 
  $\kat{S-FairCat}$.
\end{dem}

\begin{BM}
  In \ref{equiv}, equimorphisms of fair $\SS$-categories are defined,
  and in \ref{equiv-ditto} it is shown that the notion is compatible with 
  the    standard discrete objects.
\end{BM}

% \scriptsize
% 
% The existence of arbitrary sums is not strictly needed, if you 
% are willing to give up the decomposition of the space of all arrows
% into hom spaces, $A \simeq \coprod_{x,y\in O} A(x,y)$.  All we need is 
% really the individual hom spaces $A(x,y)$ in order to make sense of 
% the
% definition of fully faithful, it is not necessary to know that all
% these assemble into a sum or that this sum is $A$ again.
% So it is enough to assume that arbitrary sums of $*$ exist
% (these are the discrete objects), and that the inclusions $* \to 
% \coprod *$ admit pullbacks\ldots
% 
% \normalsize

\nocite{Leinster-0305049}
\nocite{Kock-Toen:0304}

\hyphenation{mathe-matisk}

\footnotesize

\noindent
Address:
% Departament
Dept.~%
de Matemàtiques, Universitat Autònoma de Barcelona,
08193 Bellaterra,
% (Barcelona)
{\sc Espanya}
% Spain

\noindent
E-mail: \texttt{kock@mat.uab.es}

\label{lastpage}

\end{document}